\documentclass[12pt]{article}
\usepackage{amssymb, amsmath, amsthm, amscd}
\usepackage[utf8]{inputenc}
\usepackage[all,cmtip]{xy}
\usepackage{enumitem}
\usepackage{setspace}
\usepackage{mathdots}

\usepackage[mathscr]{euscript}
\usepackage[colorlinks=true]{hyperref}
\hypersetup{colorlinks   = true}
\hypersetup{linkcolor=blue}
\usepackage{tikz}

\begin{document}

\mathchardef\mhyphen="2D
\newtheorem{The}{Theorem}[section]
\newtheorem{Lem}[The]{Lemma}
\newtheorem{Prop}[The]{Proposition}
\newtheorem{Cor}[The]{Corollary}
\newtheorem{Rem}[The]{Remark}
\newtheorem{Obs}[The]{Observation}
\newtheorem{SConj}[The]{Standard Conjecture}
\newtheorem{Titre}[The]{\!\!\!\! }
\newtheorem{Conj}[The]{Conjecture}
\newtheorem{Question}[The]{Question}
\newtheorem{Prob}[The]{Problem}
\newtheorem{Def}[The]{Definition}
\newtheorem{Not}[The]{Notation}
\newtheorem{Claim}[The]{Claim}
\newtheorem{Conc}[The]{Conclusion}
\newtheorem{Ex}[The]{Example}
\newtheorem{Fact}[The]{Fact}
\newtheorem{Formula}[The]{Formula}
\newtheorem{Formulae}[The]{Formulae}
\newtheorem{The-Def}[The]{Theorem and Definition}
\newtheorem{Prop-Def}[The]{Proposition and Definition}
\newtheorem{Lem-Def}[The]{Lemma and Definition}
\newtheorem{Cor-Def}[The]{Corollary and Definition}
\newtheorem{Conc-Def}[The]{Conclusion and Definition}
\newtheorem{Terminology}[The]{Note on terminology}
\newcommand{\C}{\mathbb{C}}
\newcommand{\R}{\mathbb{R}}
\newcommand{\N}{\mathbb{N}}
\newcommand{\Z}{\mathbb{Z}}
\newcommand{\Q}{\mathbb{Q}}
\newcommand{\Proj}{\mathbb{P}}
\newcommand{\Rc}{\mathcal{R}}
\newcommand{\Oc}{\mathcal{O}}
\newcommand{\Vc}{\mathcal{V}}
\newcommand{\Id}{\operatorname{Id}}
\newcommand{\pr}{\operatorname{pr}}
\newcommand{\rk}{\operatorname{rk}}
\newcommand{\del}{\partial}
\newcommand{\delbar}{\bar{\partial}}
\newcommand{\Cdot}{{\raisebox{-0.7ex}[0pt][0pt]{\scalebox{2.0}{$\cdot$}}}}
\newcommand\nilm{\Gamma\backslash G}
\newcommand\frg{{\mathfrak g}}
\newcommand{\fg}{\mathfrak g}
\newcommand{\Oh}{\mathcal{O}}
\newcommand{\Kur}{\operatorname{Kur}}
\newcommand\gc{\frg_\mathbb{C}}
\newcommand\hisashi[1]{{\textcolor{red}{#1}}}
\newcommand\dan[1]{{\textcolor{blue}{#1}}}
\newcommand\luis[1]{{\textcolor{orange}{#1}}}

\begin{center}

{\Large\bf Higher-Degree Holomorphic Contact Structures}

\end{center}

\begin{center}

{\large Hisashi Kasuya, Dan Popovici and Luis Ugarte}

\end{center}

\vspace{1ex}

\noindent{\small{\bf Abstract.} We introduce the classes of holomorphic $p$-contact manifolds and holomorphic $s$-symplectic manifolds that generalise the classical holomorphic contact and holomorphic symplectic structures. After observing their basic properties and exhibiting a wide range of examples, we give two types of general conceptual results involving the former class of manifolds: structure theorems and unobstructedness theorems. The latter type generalises to our context the classical Bogomolov-Tian-Todorov theorem for a type of small deformations of complex structures that generalise the small essential deformations previously introduced for the Iwasawa manifold and for Calabi-Yau page-$1$-$\partial\bar\partial$-manifolds.}

\vspace{1ex}

\section{Introduction}\label{section:introduction}

Let $X$ be an odd-dimensional compact complex manifold with $\mbox{dim}_\C X =n = 2p+1$. Recall that a {\it holomorphic contact structure} on $X$ is a smooth $(1,\,0)$-form $\eta\in C^\infty_{1,\,0}(X,\,\C)$ such that $\bar\partial\eta=0$ and $\eta\wedge(\partial\eta)^p\neq 0$ at every point of $X$. More generally, a holomorphic contact form may take values in a holomorphic line bundle $L$ over $X$, namely $\eta\in C^\infty_{1,\,0}(X,\, L)$. Intuitively, the $2p+1$ dimensions of $X$ are split into two classes: the direction of $\eta$ and the $2p$ directions of $(\partial\eta)^p$.

In this paper, we replace this splitting by the splitting $2p+1 = p + (p+1)$ by introducing and studying the following

\begin{Def}\label{Def:hol_p-contact_structure} Let $X$ be a compact complex manifold with $\mbox{dim}_\C X =n = 2p+1$.

\vspace{1ex}

\noindent $(1)$\, A {\bf holomorphic $p$-contact structure} on $X$ is a smooth $(p,\,0)$-form $\Gamma\in C^\infty_{p,\,0}(X,\,\C)$ such that \begin{eqnarray*}(a)\hspace{2ex} \bar\partial\Gamma=0 \hspace{6ex} \mbox{and} \hspace{6ex} (b)\hspace{2ex} \Gamma\wedge\partial\Gamma\neq 0 \hspace{2ex} \mbox{at every point of} \hspace{1ex} X.\end{eqnarray*}

\vspace{1ex}

\noindent $(2)$\, We say that $X$ is a {\bf holomorphic $p$-contact manifold} if there exists a  holomorphic $p$-contact structure $\Gamma$ on $X$.

\end{Def}

Holomorphic $p$-contact structures $\Gamma$ generalise the classical holomorphic contact structures $\eta$ in dimension $n = 2p+1 = 4l+3$: any such $\eta$ induces such a $\Gamma$ by taking $\Gamma = \eta\wedge(\partial\eta)^{s-1}$, where $p=2s-1$. (See (2) of Observation \ref{Obs:contact_p-contact_generalisation}.) This generalisation affords a considerably greater flexibility, as demonstrated, for example, by the existence of many compact holomorphic $p$-contact manifolds that do not admit any holomorphic contact structure (see class I of examples in Proposition \ref{Prop:higher-dim_Iwasawa}).

Holomorphic $p$-contact manifolds $(X,\,\Gamma)$ are {\it Calabi-Yau manifolds} (construed here as compact complex manifolds $X$ whose canonical bundle $K_X$ is {\it trivial}). Indeed, the holomorphic $(n,\,0)$-form $u_\Gamma := \Gamma\wedge\partial\Gamma$ is a non-vanishing global holomorphic section of $K_X$. 

Demailly considered in [Dem02] the situation where the manifold $X$ is K\"ahler (while in our case $X$ is far from K\"ahler -- see Observation \ref{Obs:non_E_1} below) and the form $\Gamma$ satisfies our condition (a) and assumes values in the dual $L^{-1}$ of a pseudo-effective holomorphic line bundle $L\rightarrow X$. The conclusion of the main result of [Dem02] (which fails without the K\"ahler assumption that is used crucially in the integrations by parts performed in the proof) was that $\Gamma$ induces a (possibly singular) foliation on $X$. This is obtained in [Dem02] as an easy consequence of the main by-product of that proof which implies that $\Gamma$ cannot satisfy our hypothesis (b) when $p$ is odd and $X$ is K\"ahler. So, our holomorphic $p$-contact manifolds are, in a sense, complementary to the manifolds studied in [Dem02].

\vspace{2ex}

In $\S$\ref{section:preliminaries}, we notice that a compact $p$-contact manifold $X$ is always of complex dimension $n\equiv 3$ {\it mod} $4$ (since $p$ must be odd) and is never K\"ahler and can never even have the far weaker property that its Fr\"olicher spectral sequence degenerates at $E_1$. However, many of the examples of such manifolds that we exhibit in $\S$\ref{section:examples} are {\it page-$1$-$\partial\bar\partial$-manifolds}. This manifold class was introduced and studied in [PSU20a], [PSU20b] and [PSU20c] and consists of those compact complex manifolds $X$ that support a Hodge theory in which the traditional Dolbeault cohomology groups ($=$ the spaces $E_1^{p,\,q}(X)$ lying on the first page of the Fr\"olicher spectral sequence) are replaced by the spaces $E_2^{p,\,q}(X)$ lying on the second page of this spectral sequence.

\vspace{1ex}

In $\S$\ref{subsection:examples_hol-s-symplectic}, we introduce another class of compact complex manifolds that we call {\bf holomorphic $s$-symplectic}. They generalise the classical notion of {\it holomorphic symplectic manifolds} by having complex dimension $n = 2s$ and carrying a smooth $(s,\,0)$-form $\Omega\in C^\infty_{s,\,0}(X,\,\C)$ such that \begin{eqnarray*}(i)\,\,\bar\partial\Omega = 0; \hspace{3ex} \mbox{and}  \hspace{3ex} (ii)\,\, \Omega\wedge\Omega \neq 0 \hspace{2ex} \mbox{at every point of}\hspace{1ex} X.\end{eqnarray*} It is easily seen that $s$ must be even, hence $n$ is a multiple of $4$, in this case. Once again, this generalisation brings a fair amount of flexibility to the classification theory of manifolds: in every complex dimension $4l\geq 8$, there are compact complex manifolds that carry holomorphic $s$-symplectic structures but do not carry any holomorphic symplectic structure. (See Example \ref{Ex:s-symplectic_non-hol-symplectic}.)

\vspace{1ex}

We use these manifolds to give two kinds of structure theorems. In $\S$\ref{subsection:examples_hol-s-symplectic}, we start off with an arbitrary complex parallelisable nilmanifold $Y$ of complex dimension $4l = 2s$, we notice that any such $Y$ has a natural {\it holomorphic $s$-symplectic structure} $\Omega\in C^\infty_{s,\,0}(Y,\,\C)$ and then we construct a class of {\it holomorphic $p$-contact manifolds} $(X,\,\Gamma)$ of complex dimension $4l+3$, each of them being equipped with a natural projection $\pi:X\longrightarrow Y$ by means of which we relate $\Gamma$ to $\Omega$. There is a special case where $X$ is the Cartesian product of $Y$ with the $3$-dimensional Iwasawa manifold $I^{(3)}$, but there are many more {\it holomorphic $p$-contact manifolds} $X$ induced in this way by a given {\it holomorphic $s$-symplectic manifold} $Y$ besides the product $Y\times I^{(3)}$. (See Theorem \ref{The:hol-s-symplectic_hol-p-contact}.)

Moreover, further examples of {\it holomorphic $p$-contact manifolds} can be obtained as a consequence of Proposition \ref{Prop:products}: the product  $X\times Y$ of any holomorphic $p$-contact manifold $X$ by any holomorphic $s$-symplectic manifold $Y$ is a holomorphic $(p+s)$-contact manifold.

In the other type of structure theorem that we give in $\S$\ref{section:structure-theorem} for arbitrary manifolds (i.e. they need not be nilmanifolds), the order of events is reversed: we start off with a {\it holomorphic $p$-contact manifold} $(X,\,\Gamma)$ of complex dimension $n = 2p+1 = 4l +3$ supposed to admit a surjective holomorphic submersion $\pi:X\longrightarrow Y$ over a $4l$-dimensional compact complex manifold $Y$ whose fibres are the leaves of a given smooth holomorphic foliation ${\cal H}\subset T^{1,\,0}X$. Under certain further assumptions, we prove in Theorem \ref{Thm:structure_Iwasawa-s-symplectic} that the base manifold $Y$ carries a {\it holomorphic $s$-symplectic structure} $\Omega$ (where $s=2l$) that is related to the given {\it holomorphic $p$-contact structure} $\Gamma$ of $X$ via the fibration map $\pi$ in a way similar to that of our first structure theorem outlined above. 

\vspace{1ex}

A key conceptual role in this paper is played by two coherent, torsion-free subsheaves ${\cal F}_\Gamma$ and ${\cal G}_\Gamma$ of the holomorphic tangent sheaf ${\cal O}(T^{1,\,0}X)$ of any holomorphic $p$-contact manifold $(X,\,\Gamma)$. They are defined in $\S$\ref{section:sheaves_F-G} as the sheaf of germs of holomorphic $(1,\,0)$-vector fields $\xi$ such that $\xi\lrcorner\Gamma = 0$, respectively $\xi\lrcorner\partial\Gamma = 0$. In the special case of a holomorphic contact manifold $(X,\,\eta)$, the corresponding sheaves ${\cal F}_\eta$ and ${\cal G}_\eta$, defined by the respective conditions $\xi\lrcorner\eta = 0$ and $\xi\lrcorner\partial\eta = 0$, are locally free, of complex co-rank $1$, respectively of complex rank $1$, ${\cal G}_\eta$ is even a {\it trivial} holomorphic line bundle in whose dual the contact form $\eta$ is a non-vanishing global holomorphic section and one has the direct sum splitting of holomorphic vector bundles: \begin{eqnarray}\label{eqn:sheaves-splitting_introd}T^{1,\,0}X = {\cal F}_\eta\oplus{\cal G}_\eta.\end{eqnarray}

However, in the general holomorphic $p$-contact case, the sheaves ${\cal F}_\Gamma$ and ${\cal G}_\Gamma$ need not be locally free, they have unpredictable ranks (there are even cases where ${\cal G}_\Gamma$ is the zero sheaf) and, although they are always in a direct sum, they need not fill out the tangent sheaf ${\cal O}(T^{1,\,0}X)$. (See Proposition \ref{Prop:F-G-sheaves_properties}.) Moreover, ${\cal F}_\Gamma$ is not integrable (this being the point of a $p$-contact structure), but ${\cal G}_\Gamma$ is (see Observation \ref{Obs:G_Gamma_integrability}). It is, thus, natural to wonder whether the above-mentioned structure theorem can be worded in terms of the leaves of the (possibly singular) holomorphic foliation induced by ${\cal G}_\Gamma$ or some other Frobenius-integrable subsheaf of ${\cal O}(T^{1,\,0}X)$ that has a direct-sum complement in the tangent sheaf.

Specifically, given a compact holomorphic $p$-contact manifold $(X,\,\Gamma)$ with $\mbox{dim}_\C X = n = 2p+1 = 4l +3$, it seems natural to wonder whether $X$ has the structure of a holomorphic fibration over a holomorphic $s$-symplectic manifold $Y$ (possibly with $s =2l$, but other values of $s$ ought to be possible) whose fibres are the leaves of ${\cal G}_\Gamma$ or of some other holomorphic foliation on $X$.

We are aware of two pieces of motivating evidence for such a result. Besides our above-mentioned Theorem \ref{The:hol-s-symplectic_hol-p-contact}, the other is Foreman's Theorem 6.2 in [For00]: an analogue in the setting of holomorphic contact manifolds of the Boothby-Wang fibration theorem in [BW58] for real contact structures. Specifically, in our notation, it stipulates that any holomorphic contact manifold $(X,\,\eta)$ of complex dimension $n=2p+1$, whose holomorphic line subbundle ${\cal G}_\eta\subset T^{1,\,0}X$ is {\it regular} in the sense of [Pal57], arises as a holomorphic fibration $\pi:X\longrightarrow Y$ over a holomorphic symplectic manifold $(Y,\,\Omega)$ such that $\pi^\star\Omega = d\eta$ on $X$. Foreman invokes a result of [Pal57] to argue that the base $Y$ of the fibration, defined as the space of maximal connected leaves of the foliation ${\cal G}_\eta$, is a manifold. The fibres are elliptic curves (since they are the leaves of the rank-one holomorphic foliation ${\cal G}_\eta$, each of them having a non-vanishing global vector field, namely $\xi_\eta$, the dual of $\eta$, so the Euler characteristic of each leaf equals $0$ by the Poincar\'e-Hopf Index Theorem and therefore, by the standard classification theorem for compact real surfaces, each leaf is homeomorphic to the torus $S^1\times S^1$) and the fibration is a principal torus bundle for which the contact form $\eta$ of $X$ is a connection form whose curvature form $d\eta$ is the pullback to $X$ of the holomorphic symplectic form $\Omega$ of $Y$.

Palais's regularity assumption (which seems to be a way of ensuring that the leaves of a foliation are compact) and his result in [Pal57] (according to which a regular foliation $\Theta$ of rank $m$ on an $n$-dimensional compact manifold $M$ has the property that the space $M/\Theta$ of maximal connected leaves, given the quotient topology, has the structure of an $(n-m)$-dimensional manifold) seem suited to the real context and to $C^\infty$ objects, less so to the complex context and to holomorphic objects. We believe that a great part of Foreman's proof of his Theorem 6.2 in [For00] can be reproved using complex-analytic arguments (see comments after the proof of Theorem \ref{Thm:structure_Iwasawa-s-symplectic}).

A natural replacement of Palais's above-mentioned regularity assumption and result in the complex setting seems to be the discussion in [EMS77] of conditions ensuring the existence of an upper bound on the volumes of the leaves in a compact manifold $X$ foliated by compact submanifolds. This is a key question since the volume boundedness near any given leaf is equivalent to the topology of the leaf space being Hausdorff near the given leaf (see [EMS77]). One sufficient condition for volume boundedness given in Theorem 1 of [EMS77] is the existence of a $d$-closed form $\omega$ on $X$ of degree equal to the dimension of the leaves such that $\omega$ has positive integral along each leaf. Such a form exists if, for example, $X$ is a compact K\"ahler manifold and the smooth foliation is holomorphic.

The volume boundedness question is also related to Lieberman's Theorem 1.1 in [Lie78] strengthening an earlier result by Bishop: a subset $(Z_s)_{s\in S}$ of the Barlet space ${\cal C}(X)$ of cycles of a compact complex manifold $X$ is relatively compact if and only if the volumes (measured with respect to an arbitrary Hermitian metric $\omega$ on $X$) of the $Z_s$'s are uniformly bounded (above). 

One of the difficulties of applying to our case the results in [EMS77] stems from the non-K\"ahlerianity of compact $p$-contact manifolds $(X,\,\Gamma)$ and from the arbitrary and unpredictable dimensions of the leaves of the holomorphic foliations (e.g. ${\cal G}_\Gamma$) involved. For example, when the leaves are of complex codimension $1$, one can apply the well-known result (see e.g. [CP94, Remark 2.18.]) according to which the irreducible components of the Barlet space ${\cal C}^{n-1}(X)$ of divisors of an $n$-dimensional compact complex manifold $X$ are {\bf compact}. However, in our case, if our numerous examples are anything to go by, the leaves are almost never of codimension $1$.

Thus, the general question of when the leaves of a holomorphic foliation are the fibres of a holomorphic fibration, despite its having been studied for several decades, is still quite mysterious outside a few special cases and deserves a separate treatment. For this reason, in our case, we will suppose the existence of a fibration whose fibres are the leaves of a given foliation and will concentrate on finding sufficient conditions for a given holomorphic $p$-contact structure on the total space giving rise to a holomorphic $s$-symplectic structure on the base of the fibration.

\vspace{1ex}

In fact, we obtain two types of general conceptual results in relation to our manifolds introduced in this paper. One consists in the structure theorems mentioned above, given in $\S$\ref{subsection:examples_hol-s-symplectic} and in $\S$\ref{section:structure-theorem}.

Another consists of general conceptual results pertaining to small deformations of holomorphic $p$-contact manifolds $(X,\,\Gamma)$ and occupying $\S$\ref{section:contact-def}.

A non-K\"ahler mirror symmetry project was initiated in [Pop18] starting from the introduction of {\it small essential deformations} of the Iwasawa manifold $I^{(3)}$. These form a subclass, parametrised by the space $E_2^{n-1,\,1}$ lying on the second page of the Fr\"olicher spectral sequence, of the Kuranishi family of $I^{(3)}$ (that is classically known to be parametrised by the space $E_1^{n-1,\,1}$ lying on the first page). In that case, $n= \mbox{dim}_\C I^{(3)} =3$. The essential small deformations were shown to correspond, via an explicitly defined mirror map that was then proved to be a local biholomorphism, to the complexified Gauduchon cone of $I^{(3)}$.

This approach to mirror symmetry for certain Calabi-Yau manifolds in the general, possibly non-K\"ahler, setting was then pursued in [PSU20c] through the extension of the notion of {\it small essential deformations} to the class of Calabi-Yau {\it page-$1$-$\partial\bar\partial$-manifolds} that had been introduced in [PSU20a] and to which $I^{(3)}$ belongs.

We now further pursue this strategy in $\S$\ref{section:contact-def} by defining and studying an analogue of the small essential deformations in the context of holomorphic $p$-contact manifolds $(X,\,\Gamma)$. In Definition \ref{Def:essential-horizontal-def}, we propose the notion of {\bf essential horizontal deformations} of $(X,\,\Gamma)$ which we then prove in Theorem \ref{The:full-unobstructedness} to be {\bf unobstructed} in the usual sense adapted to our context.


Our unobstructedness result generalises in the context of holomorphic $p$-contact manifolds the classical Bogomolov-Tian-Todorov theorem stipulating the unobstructedness of the small deformations of any compact K\"ahler (or merely $\partial\bar\partial$) Calabi-Yau manifold.

\vspace{1ex}

This discussion throws up some enticing ideas, one of which being the prospect of building a mirror symmetry theory for compact holomorphic $p$-contact manifolds (possibly assumed to satisfy some extra conditions, such as having the {\it page-$1$-$\partial\bar\partial$-property}) that would exchange the roles of $\Gamma$ and $\partial\Gamma$ (or those of ${\cal F}_\Gamma$ and ${\cal G}_\Gamma$, or those of some other pair of subsheaves or subbundles of the holomorphic tangent bundle) when switching from one member to the other of any pair $(X,\,\widetilde{X})$ of mirror symmetric holomorphic $p$-contact manifolds. For example, what would become of the structure Theorem \ref{Thm:structure_Iwasawa-s-symplectic} or of the small {\it  essential horizontal deformations} when crossing the mirror from $(X,\,\Gamma)$ to $(\widetilde{X},\,\widetilde{\Gamma})$?

\section{Preliminaries}\label{section:preliminaries} We will make two kinds of observations about these classes of manifolds.

\subsection{Cohomology}\label{subsection:cohomology}

We first notice that holomorphic $p$-contact manifolds lie a long way away from the K\"ahler realm. 

\begin{Obs}\label{Obs:non_E_1} The Fr\"olicher spectral sequence of no holomorphic $p$-contact manifold $X$ degenerates at $E_1$.  

  In particular, no $\partial\bar\partial$-manifold (hence no compact K\"ahler manifold) carries a holomorphic $p$-contact structure.

\end{Obs}  

\noindent {\it Proof.} Suppose there exists a holomorphic $p$-contact structure $\Gamma$ on a compact $(2p+1)$-dimensional complex manifold $X$. Then, $\Gamma$ represents a Dolbeault cohomology class $[\Gamma]_{\bar\partial}\in H^{p,\,0}_{\bar\partial}(X,\,\C)\simeq E_1^{p,\,0}(X)$ whose image under the differential \begin{eqnarray*}d_1:E_1^{p,\,0}(X)\longrightarrow E_1^{p+1,\,0}(X)\end{eqnarray*} on the first page of the  Fr\"olicher spectral sequence of $X$ is $d_1([\Gamma]_{\bar\partial}) = [\partial\Gamma]_{\bar\partial}\in E_1^{p+1,\,0}(X)$. 

If the Fr\"olicher spectral sequence degenerated at $E_1$, all the differentials $d_1:E_1^{r,\,s}(X)\longrightarrow E_1^{r+1,\,s}(X)$ on its first page would be identically zero. In particular, we would have $[\partial\Gamma]_{\bar\partial} = 0\in E_1^{p+1,\,0}(X)$, so there would exist $\zeta\in C^\infty_{p+1,\,-1}(X,\,\C)$ such that $\partial\Gamma = \bar\partial\zeta$. However, for bidegree reasons, we would then have $\zeta = 0$, hence $\partial\Gamma = 0$, which is ruled out by property (b) under $(1)$ of Definition \ref{Def:hol_p-contact_structure}. \hfill $\Box$

\vspace{2ex}

 Thus, the best possible degeneration of the Fr\"olicher spectral sequence of any holomorphic $p$-contact manifold $X$ can occur at the second page.

We will sometimes make the stronger assumption that $X$ be a {\it page-$1$-$\partial\bar\partial$-manifold}. This property was introduced in [PSU20a] and means that the Fr\"olicher spectral sequence of $X$ degenerates at the second page (a property denoted by $E_2(X) = E_\infty(X)$) and that the De Rham cohomology of $X$ is {\it pure}. Equivalently, it means that $X$ has a canonical Hodge decomposition in which the $E_1$-spaces $E_1^{p,\,q}(X)$ ($=$ the Dolbeault cohomology groups $H^{p,\,q}(X,\,\C)$ of $X$) have been replaced by the $E_2$-spaces $E_2^{p,\,q}(X)$ in the sense that the identity map induces {\it canonical isomorphisms} \begin{eqnarray}\label{eqn:Hodge-decomp_E_2}H^k_{dR}(X,\,\C)\longrightarrow\bigoplus\limits_{p+q=K}E_2^{p,\,q}(X), \hspace{5ex} k=0,1,\dots , 2n.\end{eqnarray}

Explicitly, $X$ being a {\it page-$1$-$\partial\bar\partial$-manifold} means that each cohomology class in every space $E_2^{p,\,q}(X)$ can be represented by a smooth $d$-closed $(p,\,q)$-form and that, for every $k\in\{0,\dots , 2n\}$, the linear map \begin{eqnarray*}\bigoplus\limits_{p+q=k}E_2^{p,\,q}(X)\ni\bigg(\{\alpha^{p,\,q}\}_{E_2}\bigg)_{p+q=k}\longmapsto\bigg\{\sum\limits_{p+q=k}\alpha^{p,\,q}\bigg\}_{dR}\in H^k_{dR}(X,\,\C)\end{eqnarray*} is {\it well defined} (in the sense that it is independent of the choices of $d$-closed pure-type representatives $\alpha^{p,\,q}$ of their respective $E_2$-classes) and {\it bijective}.

Recall that the description of the Fr\"olicher spectral sequence of any compact complex manifold $X$ given in [CFGU97] reinterprets the $\C$-vector spaces $E_2^{p,\,q}(X)$ featuring on the second page as the quotients \begin{eqnarray*}E_2^{p,\,q}(X) = \frac{{\cal Z}_2^{p,\,q}(X)}{{\cal C}_2^{p,\,q}(X)}, \hspace{6ex} p,\,q=0,\dots , n,\end{eqnarray*} where the $\C$-vector space ${\cal Z}_2^{p,\,q}(X)$ consists of the $C^\infty$ $(p,\,q)$-forms $u$ such that \begin{eqnarray}\label{eqn:Z_2_def}\bar\partial u=0 \hspace{5ex} \mbox{and} \hspace{5ex} \partial u\in\mbox{Im}\,\bar\partial\end{eqnarray} (these forms $u$ being called {\it $E_2$-closed}), while the $\C$-vector space ${\cal C}_2^{p,\,q}(X)$ consists of the $C^\infty$ $(p,\,q)$-forms $u$ for which there exist forms $\zeta\in C^\infty_{p-1,\,q}(X,\,\C)$ and $\xi\in C^\infty_{p,\,q-1}(X,\,\C)$ such that \begin{eqnarray*}u = \partial\zeta + \bar\partial\xi \hspace{5ex} \mbox{and} \hspace{5ex} \bar\partial\zeta = 0\end{eqnarray*} (these forms $u$ being called {\it $E_2$-exact}). We always have the inclusions: \begin{eqnarray*}{\cal C}_1^{p,\,q}(X)\subset{\cal C}_2^{p,\,q}(X)\subset{\cal Z}_2^{p,\,q}(X)\subset{\cal Z}_1^{p,\,q}(X),\end{eqnarray*} where ${\cal C}_1^{p,\,q}(X)$ (resp. ${\cal Z}_1^{p,\,q}(X)$) is the $\C$-vector space of $C^\infty$ $(p,\,q)$-forms $u$ such that $u\in\mbox{Im}\,\bar\partial$ (resp. $u\in\ker\bar\partial$). Obviously, $E_1^{p,\,q}(X) = {\cal Z}_1^{p,\,q}(X)/{\cal C}_1^{p,\,q}(X)$ for all $p,\,q$.

In this paper, we will use the following consequence of Theorem 4.3. in [PSU20b]: if $X$ is a {\it page-$1$-$\partial\bar\partial$-manifold}, the following equality of $\C$-vector spaces holds in every bidegree $(p,\,q)$: \begin{eqnarray}\label{eqn:ddbar_im_Z-2}\partial\bigg({\cal Z}_2^{p,\,q}(X)\bigg) = \mbox{Im}\,(\partial\bar\partial).\end{eqnarray} The inclusion ``$\supset$'' holds on any $X$, it is the inclusion ``$\subset$'' that follows from the {\it page-$1$-$\partial\bar\partial$}-assumption.

\subsection{Dimension}\label{subsection:dimension} The next observation is that the complex dimensions of compact holomorphic $p$-contact manifolds are congruent to $3$ modulo $4$.   

\begin{Obs}\label{Obs:dim_4l-plus-3} Let $X$ be a compact complex manifold with $\mbox{dim}_\C X = n = 2p+1$. If there exists a {\bf holomorphic $p$-contact structure} $\Gamma\in C^\infty_{p,\,0}(X,\,\C)$ on $X$, $p$ is {\bf odd}.

\end{Obs}  

\noindent {\it Proof.} If $p$ were even, we would get: $\frac{1}{2}\,\partial(\Gamma^2) = \Gamma\wedge\partial\Gamma$. This is a non-vanishing holomorphic $(n,\,0)$-form on $X$ since $\Gamma$ is a holomorphic $p$-contact structure. Hence \begin{eqnarray*}dV_\Gamma := i^{n^2}\,(\Gamma\wedge\partial\Gamma)\wedge(\overline\Gamma\wedge\bar\partial\overline\Gamma) = \frac{i^{n^2}}{4}\,\partial\bigg(\Gamma^2\wedge\bar\partial\overline\Gamma^2\bigg)\end{eqnarray*} would be a strictly positive $\partial$-exact $(n,\,n)$-form on $X$.

We would then get (using Stokes to infer the equality below): \begin{eqnarray*}0<\int_X dV_\Gamma = 0,\end{eqnarray*} a contradiction.  \hfill $\Box$

\subsection{Calabi-Yau property}\label{subsection:C-Y} Let $X$ be a compact complex manifold with $\mbox{dim}_\C X = n=2p+1$. Suppose there exists a {\it $p$-contact structure} $\Gamma$ on $X$.

Since $\Gamma\wedge\partial\Gamma$ is a non-vanishing holomorphic $(n,\,0)$-form on $X$, it defines a non-vanishing global holomorphic section of the canonical line bundle $K_X$. Consequently, $K_X$ is {\it trivial}, so $X$ is a Calabi-Yau manifold. Moreover, the smooth $(n,\,n)$-form \begin{eqnarray*}dV_\Gamma:=(i^{p^2}\,\Gamma\wedge\overline\Gamma)\wedge(i^{(p+1)^2}\,\partial\Gamma\wedge\bar\partial\overline\Gamma) = i^{n^2}\,(\Gamma\wedge\partial\Gamma)\wedge(\overline\Gamma\wedge\bar\partial\overline\Gamma)\end{eqnarray*} is strictly positive on $X$, so it defines a volume form induced by the contact structure $\Gamma$. After possibly multiplying $\Gamma$ by a constant, we may assume that \begin{eqnarray*}\int\limits_X dV_\Gamma =1.\end{eqnarray*} After this normalisation has been performed, we will call the non-vanishing holomorphic $(n,\,0)$-form $u_\Gamma:=\Gamma\wedge\partial\Gamma$ the {\it Calabi-Yau form} induced by $\Gamma$.

\vspace{2ex}

Let $\bar\partial$ still denote the canonical $\bar\partial$-operator defining the holomorphic structure of the holomorphic tangent vector bundle $T^{1,\,0}X$. Meanwhile, the contraction operation $\xi\lrcorner\cdot\,:C^\infty_{p,\,q}(X,\,\C)\longrightarrow C^\infty_{p-1,\,q}(X,\,\C)$ by any vector field $\xi\in C^\infty(X,\,T^{1,\,0}X)$ applies pointwise to $\C$-valued forms of any bidegree $(p,\,q)$. The general formula \begin{eqnarray}\label{eqn:xi_contraction-Leibniz}\bar\partial(\xi\lrcorner u) = (\bar\partial\xi)\lrcorner u - \xi\lrcorner(\bar\partial u)\end{eqnarray} holds for any $u\in C^\infty_{p,\,q}(X,\,\C)$ -- see e.g. [Pop19, Lemma 4.3].

Similarly, for any $\theta\in C^\infty_{0,\,1}(X,\,T^{1,\,0}X)$, the pointwise operation $\theta\lrcorner\cdot\,:C^\infty_{p,\,q}(X,\,\C)\longrightarrow C^\infty_{p-1,\,q+1}(X,\,\C)$ combines the contraction of any form $u\in C^\infty_{p,\,q}(X,\,\C)$ by the $T^{1,\,0}X$-part of $\theta$ with the multiplication by the $(0,\,1)$-form part of $\theta$. The general formula \begin{eqnarray}\label{eqn:theta_contraction-Leibniz}\bar\partial(\theta\lrcorner u) = (\bar\partial\theta)\lrcorner u + \theta\lrcorner(\bar\partial u)\end{eqnarray} holds for any $u\in C^\infty_{p,\,q}(X,\,\C)$ -- see e.g. [Pop19, Lemma 4.3].

\vspace{2ex}

Finally, recall the {\it Calabi-Yau isomorphism} (defined for any non-vanishing holomorphic $(n,\,0)$-form $u$, in particular for our $u_\Gamma$, and for every $q\in\{0,1,\dots ,n\}$): \begin{eqnarray}\label{eqn:C_Y-isomorphism_scalar}T_\Gamma : C^\infty_{0,\,q}(X,\,T^{1,\,0}X)\stackrel{\,\cdot\lrcorner u_\Gamma\hspace{2ex}}{\longrightarrow} C^\infty_{n-1,\,q}(X,\,\C),  \hspace{3ex} T_\Gamma(\theta):= \theta\lrcorner u_\Gamma = \theta\lrcorner(\Gamma\wedge\partial\Gamma)\end{eqnarray} and the isomorphism it induces in cohomology \begin{eqnarray}\label{eqn:C_Y-isomorphism_cohom}T_{[\Gamma]} : H^{0,\,q}_{\bar\partial}(X,\,T^{1,\,0}X)\stackrel{\,\cdot\lrcorner[u_\Gamma]\hspace{2ex}}{\longrightarrow} H^{n-1,\,q}_{\bar\partial}(X,\,\C),  \hspace{3ex} T_{[\Gamma]}([\theta]_{\bar\partial}):= [\theta\lrcorner u_\Gamma]_{\bar\partial} = [\theta\lrcorner(\Gamma\wedge\partial\Gamma)]_{\bar\partial}.\end{eqnarray} Note that our hypotheses and formulae (\ref{eqn:xi_contraction-Leibniz})--(\ref{eqn:theta_contraction-Leibniz}) imply the well-definedness of $T_{[\Gamma]}$ (i.e. the fact that all the cohomology classes featuring in (\ref{eqn:C_Y-isomorphism_cohom}) are independent of the choice of representative $\theta$ of $[\theta]_{\bar\partial}$.)

\subsection{Lie derivatives with respect to vector-valued forms}\label{subsection:Lie}

Most of the results in this subsection are proved in $\S4.3$ and $\S4.4$ of [KPU25]. We include the statements here for the reader's convenience.  

\vspace{1ex}

Recall that the Lie bracket between two elements $\varphi\in C^\infty_{0,\,p}(X,\,T^{1,\,0}X)$ and $\psi\in C^\infty_{0,\,q}(X,\,T^{1,\,0}X)$, whose expressions in local holomorphic coordinates are \begin{eqnarray*}\varphi = \sum\limits_{\lambda=1}^n\varphi^\lambda\,\frac{\partial}{\partial z_\lambda}  \hspace{3ex} \mbox{and} \hspace{3ex} \psi = \sum\limits_{\lambda=1}^n\psi^\lambda\,\frac{\partial}{\partial z_\lambda},\end{eqnarray*} where the $\varphi^\lambda$'s, resp. the $\psi^\lambda$'s, are $\C$-valued $(0,\,p)$-forms, resp. $\C$-valued $(0,\,q)$-forms, is defined by \begin{eqnarray}\label{eqn:bracket_0pq-vector_def}[\varphi,\,\psi]:=\sum\limits_{\lambda,\,\mu=1}^n\bigg(\varphi^\mu\wedge\frac{\partial\psi^\lambda}{\partial z_\mu} - (-1)^{pq}\,\psi^\mu\wedge\frac{\partial\varphi^\lambda}{\partial z_\mu}\bigg)\,\frac{\partial}{\partial z_\lambda}.\end{eqnarray} 
 Thus, $[\varphi,\,\psi]\in C^\infty_{0,\,p+q}(X,\,T^{1,\,0}X)$. 
 For the basic properties of this Lie bracket see e.g. [Kod86, $\S5$]. 
 
 \vspace{1ex}

With every $\theta\in C^{\infty}_{0,\,q}(X,\,T^{1,\,0}X)$, we associate the linear differential operator $L_\theta:C^\infty_{r,\,s}(X,\,\C)\longrightarrow C^\infty_{r,\,s+q}(X,\,\C)$ of order one and bidegree $(0,\,q)$ defined so that the following analogue of the Leibniz formula holds in this case: \begin{equation*}\partial(\theta\lrcorner u) = L_\theta u + (-1)^{q-1}\,\theta\lrcorner\partial u.\end{equation*} The above sign is due to the total degree of any such $\theta$ being $q-1$, with the $T^{1,\,0}X$-part counting for $-1$. The case $q=0$ was introduced and used in [PU23], where the Lie derivative operator with respect to a $(1,\,0)$-vector field $\xi\in C^{\infty}(X,\,T^{1,\,0}X)$ was denoted by $L^{1,\,0}_{\xi} :C^\infty_{r,\,s}(X,\,\C)\longrightarrow C^\infty_{r,\,s}(X,\,\C)$.

Among the properties observed in Lemma 4.2 of [PU23] for the Lie derivative $L^{1,\,0}_\xi$ are the following identities: \begin{eqnarray}\label{eqn:L_commutations} (a)\,\, [\xi\lrcorner\cdot,\,L^{1,\,0}_{\eta}] = [L^{1,\,0}_{\xi},\,\eta\lrcorner\cdot] = [\xi,\,\eta]\lrcorner\cdot \hspace{2ex} & \mbox{and} & \hspace{2ex} (b)\,\,[L^{1,\,0}_{\xi},\,L^{1,\,0}_{\eta}] = L^{1,\,0}_{[\xi,\,\eta]} \\
 \nonumber   (c)\,\,[L^{1,\,0}_{\xi},\,\bar\partial]= [\partial,\,\bar\partial\xi\lrcorner\cdot\,] \hspace{2ex} & \mbox{and} & \hspace{2ex} (d)\,\, L^{1,\,0}_{\xi}(u\wedge v) = (L^{1,\,0}_{\xi}u)\wedge v + u\wedge L^{1,\,0}_{\xi}v  \end{eqnarray} that hold for all $(1,\,0)$-vector fields $\xi,\eta\in C^{\infty}(X,\,T^{1,\,0}X)$ and all $\C$-valued differential forms $u,v$ (of any degrees). In particular, if $\xi$ is holomorphic, (c) becomes $[L^{1,\,0}_{\xi},\,\bar\partial] = 0$.

\begin{Lem}\label{Lem:Lie-derivatives-mixed_prop} Let $\theta,\psi\in C^{\infty}_{0,\,1}(X,\,T^{1,\,0}X)$. The following identities hold: \begin{eqnarray}\label{eqn:Lie-derivatives-mixed_prop_0} & & [L_\theta,\,\bar\partial] = - [(\bar\partial\theta)\lrcorner\cdot\,,\,\partial]; \hspace{2ex} \mbox{in particular},\hspace{2ex} L_\theta\bar\partial = -\bar\partial L_\theta \hspace{2ex} \mbox{when} \hspace{2ex} \bar\partial\theta = 0; \\
    \label{eqn:Lie-derivatives-mixed_prop_0-bis} & & L_\theta\partial = -\partial L_\theta; \\
    \label{eqn:Lie-derivatives-mixed_prop_1} & & L_\theta(u\wedge v) = L_\theta(u)\wedge v + (-1)^{\deg u}\,u\wedge L_\theta(v) \hspace{5ex} \mbox{for all forms}\hspace{1ex} u,v;  \\
    \label{eqn:Lie-derivatives-mixed_prop_2} & & [\theta\lrcorner\cdot\,,\,L_\psi] = -[L_\theta,\,\psi\lrcorner\cdot\,] = [\theta,\,\psi]\lrcorner\cdot\, ; \\
   \label{eqn:Lie-derivatives-mixed_prop_3} & & [L_\theta\,,\,L_\psi] = L_{[\theta,\,\psi]}.\end{eqnarray}

\vspace{1ex}

\end{Lem}

\noindent {\it Proof.} See [KPU25, Lemma 4.4.]. \hfill $\Box$

\vspace{2ex}

The following terminology will come in handy later on.

\begin{Prop-Def}\label{Def:constantly_horizontal-vertical} Let $X$ be a complex manifold with $\mbox{dim}_\C X = n = 2p+1$. Suppose that $X$ carries a holomorphic $p$-contact structure $\Gamma$. Let $U\subset X$ be an open subset and let $q\in\{0,\dots , n\}$.

\vspace{1ex}

$(1)$\, A form $\theta\in C_{0,\,q}^\infty(U,\,T^{1,\,0}X)$ is said to be {\bf constantly horizontal} if \begin{eqnarray*}(a)\,\theta\lrcorner\Gamma = 0 \hspace{3ex}\mbox{and}\hspace{3ex} (b)\, L_\theta(\partial\Gamma) = 0.\end{eqnarray*} 

These conditions are equivalent to $\theta\lrcorner\Gamma = 0$ and $\partial(\theta\lrcorner\partial\Gamma) = 0$.

\vspace{1ex}

$(2)$\, A form $\theta\in C_{0,\,q}^\infty(U,\,T^{1,\,0}X)$, is said to be {\bf constantly vertical} if \begin{eqnarray*}(c)\,\theta\lrcorner\partial\Gamma = 0 \hspace{3ex}\mbox{and}\hspace{3ex} (d)\, L_\theta(\Gamma) = 0.\end{eqnarray*}

These conditions are equivalent to $\theta\lrcorner\partial\Gamma = 0$ and $\partial(\theta\lrcorner\Gamma) = 0$.

\end{Prop-Def}  

\noindent {\it Proof.} See [KPU25, Proposition and De???nition 4.7.]. \hfill $\Box$

\vspace{2ex}

By way of a justification of the terminology introduced in the above definition, we note that condition (a) in $(1)$ requires $\theta$ to be {\it horizontal} (i.e. $\theta$ to be an ${\cal F}_\Gamma$-valued $(0,\,q)$-form), while condition (b) requires $\partial\Gamma$ to be {\it constant} in the $\theta$-direction. The roles of $\Gamma$ and $\partial\Gamma$ get permuted and the term {\it vertical} is substituted for the term {\it horizontal} (or, equivalently, the sheaf ${\cal G}_\Gamma$ is substituted for the sheaf ${\cal F}_\Gamma$) in passing from $(1)$ to $(2)$.

\vspace{2ex}

With this terminology in place, we get the following

\begin{Cor}\label{Cor:brackets_horizontal-vertical_forms} Let $X$ be a complex manifold with $\mbox{dim}_\C X = n =2p+1$. Suppose that $X$ carries a holomorphic $p$-contact structure $\Gamma\in C^\infty_{p,\,0}(X,\,\C)$. Let $\theta_0,\theta_1\in C^\infty_{0,\,1}(X,\,T^{1,\,0}X)$.

\vspace{1ex}

(i)\, If one of the forms $\theta_0$ and $\theta_1$ is {\bf constantly horizontal} and the other {\bf vertical}, then $[\theta_0,\,\theta_1]$ is {\bf vertical}.

Similarly, if one of the forms $\theta_0$ and $\theta_1$ is {\bf constantly vertical} and the other {\bf horizontal}, then $[\theta_0,\,\theta_1]$ is {\bf horizontal}.

Finally, if one of the forms $\theta_0$ and $\theta_1$ is {\bf constantly horizontal} and the other {\bf constantly vertical}, then $[\theta_0,\,\theta_1] = 0$.
  
\vspace{1ex}

(ii)\, If $\theta_0$ and $\theta_1$ are both {\bf constantly horizontal} and if either $L_{\theta_1}(\theta_0\lrcorner\partial\Gamma) = 0$ or $L_{\theta_0}(\theta_1\lrcorner\partial\Gamma) = 0$, then $[\theta_0,\,\theta_1]$ is {\bf vertical}.

\end{Cor}

\noindent {\it Proof.} See [KPU25, Corollary 4.8.]. \hfill $\Box$

\vspace{2ex}

The properties of the Lie derivative w.r.t. a $T^{1,\,0}X$-valued $(0,\,1)$-form introduced and studied in [KPU25] yield at once a more general form of the classical Tian-Todorov lemma 
(cf. Lemma 3.1. in [Tia87], Lemma 1.2.4. in [Tod89])
in which the scalar-valued forms $\theta_1\lrcorner u$ and $\theta_2\lrcorner u$ are not assumed $\partial$-closed. We need not even suppose $K_X$ to be trivial.

\begin{Lem}(generalised Tian-Todorov)\label{Lem:generalised_Tian-Todorov} Let $X$ be a compact complex manifold with $n=\mbox{dim}_{\C}X$.

\vspace{1ex}

 For any bidegree $(p,\,q)$, any form $\alpha\in C^\infty_{p,\,q}(X,\,\C)$ and any forms $\theta_1, \theta_2\in C^{\infty}_{0,\, 1}(X,\, T^{1,\, 0}X)$, the following identity holds: \begin{equation}\label{eqn:generalised_basic-trick_non-C-Y}[\theta_1,\, \theta_2]\lrcorner\alpha = -\partial\bigg(\theta_1\lrcorner(\theta_2\lrcorner\alpha)\bigg) + \theta_1\lrcorner L_{\theta_2}(\alpha) + \theta_2\lrcorner L_{\theta_1}(\alpha) + \theta_1\lrcorner(\theta_2\lrcorner\partial\alpha).\end{equation}

In particular, if $\partial\alpha = 0$ and $\partial(\theta_1\lrcorner\alpha) = \partial(\theta_2\lrcorner\alpha) = 0$, then \begin{equation}\label{eqn:generalised_basic-trick_non-C-Y_re}[\theta_1,\, \theta_2]\lrcorner\alpha = -\partial\bigg(\theta_1\lrcorner(\theta_2\lrcorner\alpha)\bigg).\end{equation}

\end{Lem}

\noindent {\it Proof.} See [KPU25, (ii) of Lemma 4.11.]. \hfill $\Box$

 \section{Examples}\label{section:examples} The first observation, that implicitly provides examples, is that our holomorphic $p$-contact structures generalise the standard holomorphic contact structures in all the dimensions $n\equiv 3$ {\it mod} $4$.

 \begin{Obs}\label{Obs:contact_p-contact_generalisation} Let $X$ be a compact complex manifold with $\mbox{dim}_\C X = n = 2p+1$ and $p$ odd.

\vspace{1ex}

$(1)$\, If $n=3$ (and $p=1$), a form $\Gamma\in C^\infty_{1,\,0}(X,\,\C)$ is a holomorphic contact structure on $X$ if and only if it is a holomorphic $1$-contact structure.

\vspace{1ex}

$(2)$\, Let $s$ be the positive integer such that $p=2s-1$. Then, for any holomorphic contact structure $\eta\in C^\infty_{1,\,0}(X,\,\C)$ on $X$, the form \begin{eqnarray*}\Gamma:=\eta\wedge(\partial\eta)^{s-1}\in C^\infty_{p,\,0}(X,\,\C)\end{eqnarray*} is a holomorphic $p$-contact structure on $X$.

\end{Obs}   

 \noindent {\it Proof.} $(1)$ follows at once from the definitions. To prove $(2)$, we notice that $\bar\partial\eta = 0$ implies $\bar\partial\Gamma = 0$, while $\partial\Gamma = (\partial\eta)^s$, so $\Gamma\wedge\partial\Gamma = \eta\wedge(\partial\eta)^p\neq 0$ at every point of $X$ since $\eta$ is a contact structure.  \hfill $\Box$

\subsection{Examples in arbitrary dimensions $n\equiv 3$ {\it mod} $4$}\label{subsection:examples_dimension_4l+3} We now present two classes of examples consisting of higher-dimensional generalisations of the classical $3$-dimensional Iwasawa manifold $I^{(3)}$. All these Iwasawa-type manifolds that we discuss here have dimensions $n\equiv 3$ {\it mod} $4$ and are obtained by modifying the {\it page}-$1$-$\partial\bar\partial$-{\it manifolds}, that are also nilmanifolds endowed with invariant {\it abelian} complex structures, presented in [PSU20a, Theorem 4.8]. The modification consists in suppressing the conjugates on the r.h.s. of the structure equations, a fact that destroys the {\it abelian} property of the invariant complex structure. 

We note that the first class of higher dimensional holomorphic $p$-contact examples cannot be derived from a holomorphic contact structure on the manifold.

 \begin{Prop}\label{Prop:higher-dim_Iwasawa} Let $n = 2p+1 = 4l + 3$, with $p$ and $l$ positive integers, and let $G$ be the nilpotent $n$-dimensional complex Lie group equipped with the complex structure defined by either of the following two classes of structure equations involving a basis of holomorphic $(1,\,0)$-forms $\varphi_1,\dots , \varphi_n$:

   \vspace{1ex}

   {\bf (Class I)} \hspace{5ex} $d\varphi_1 = d\varphi_2 = 0, \hspace{2ex} d\varphi_3 = \varphi_2\wedge\varphi_1, \dots ,  d\varphi_n = \varphi_{n-1}\wedge\varphi_1$;

   \vspace{1ex}

   {\bf (Class II)} \hspace{5ex} $d\varphi_1 = \dots = d\varphi_{n-1} = 0, \hspace{2ex} d\varphi_n = \varphi_1\wedge\varphi_2 + \varphi_3\wedge\varphi_4 + \dots + \varphi_{n-2}\wedge\varphi_{n-1}$.

   \vspace{2ex}

   Then:

   \vspace{1ex}
   
$(1)$\, if $G$ is in {\bf class I}, given any $(4l+1)\times(4l+1)$ invertible upper triangular matrix $A=(a_{ij})_{3\leq i,j\leq 4l+3}$, we define the $(1,\,0)$-forms $\gamma_u= \sum_{i=u}^{4l+3} a_{ui}\,\varphi_i$, for every $u\geq 3$, and the $(p=2l+1,\,0)$-form $\Gamma_{l}$ given by
$$
\Gamma_{l} = \gamma_3\wedge(\gamma_4\wedge\gamma_5+\gamma_5\wedge\gamma_6+\gamma_6\wedge\gamma_7)\wedge\dots\wedge(\gamma_{4l}\wedge\gamma_{4l+1}+\gamma_{4l+1}\wedge\gamma_{4l+2}+\gamma_{4l+2}\wedge\gamma_{4l+3}).
$$
For any co-compact lattice $\Lambda\subset G$, the $(p=2l+1,\,0)$-form $\Gamma_{l}$ defines a {\bf holomorphic $p$-contact structure} on the compact nilmanifold $X=G/\Lambda$. 
Furthermore, $X$ does not admit any holomorphic contact structure.

   \vspace{1ex}
   
   $(2)$\, if $G$ is in {\bf class II}, for any co-compact lattice $\Lambda\subset G$, the $(p,\,0)$-form \begin{eqnarray*}\Gamma = \varphi_n\wedge(\partial\varphi_n)^l = \varphi_n\wedge\bigg(\sum\limits_{j=0}^{2l}\varphi_{2j+1}\wedge\varphi_{2j+2}\bigg)^l\end{eqnarray*} defines a {\bf holomorphic $p$-contact structure} on the compact nilmanifold $X=G/\Lambda$. This is induced by the {\bf standard contact structure} $\varphi_n$ on $X$.

\end{Prop}   

\noindent {\it Proof.} $(2)$\, Since $\bar\partial\varphi_j = 0$ for every $j$, $\bar\partial\Gamma = 0$. Moreover, $\partial\Gamma = (\partial\varphi_n)^{l+1}$, hence $$\Gamma\wedge\partial\Gamma = \varphi_n\wedge(\partial\varphi_n)^p = \varphi_n\wedge\bigg(\sum\limits_{j=0}^{2l}\varphi_{2j+1}\wedge\varphi_{2j+2}\bigg)^p$$ and this form is non-vanishing at every point of $X$ since $(\partial\varphi_n)^p$ is a sum of products of forms chosen only from $\varphi_1,\dots , \varphi_{n-1}$.

\vspace{1ex}

$(1)$\, Since $X$ is a complex parallelisable nilmanifold, any holomorphic $1$-form $\gamma$ is left invariant, so we can write $\gamma=\sum_{i=1}^{4l+3} \lambda_i\,\varphi_i$, for some $\lambda_i\in\C$. Now, $\partial \gamma = \big(\sum_{i=3}^{4l+3} \lambda_{i}\,\varphi_{i-1}\big)\wedge \varphi_1$, which implies $(\partial \gamma)^2=0$. Since $l\geq 1$, we have $p\geq 2$ and $\gamma\wedge (\partial \gamma)^p=0$, so $X$ does not admit any holomorphic contact structure.

Since all the forms $\gamma_u$ are holomorphic, $\Gamma_l$ too is holomorphic, so we have to prove that $\Gamma_l\wedge \partial\Gamma_l\not=0$. We will prove the result by induction on $l\geq 1$. Let us consider the $(1,\,0)$-forms 
$\eta_u^v= \sum_{i=u}^v a_{ui}\,\varphi_i$ for every $u\geq 3$ and $v\geq u$ and the $(2l+1,\,0)$-forms $\Phi_{l}^v$, with $v\geq 4l+3$, given by 
$$
\Phi_{l}^v = \eta_3^{v}\wedge \big(\eta_4^{v}\wedge\eta_5^{v}+\eta_5^{v}\wedge\eta_6^{v}+\eta_6^{v}\wedge\eta_7^{v}\big) \wedge\dots\wedge \big(\eta_{4l}^{v}\wedge\eta_{4l+1}^{v}+\eta_{4l+1}^{v}\wedge\eta_{4l+2}^{v}+\eta_{4l+2}^{v}\wedge\eta_{4l+3}^{v}\big).
$$
Note that $\gamma_u=\eta_u^{4l+3}$, so $\Gamma_l=\Phi_l^{4l+3}$. 
For $l=1$, we have  
\begin{eqnarray*}
\Gamma_1\wedge\partial \Gamma_1 & = & \eta_3^{7}\wedge\partial\eta_3^{7}\wedge\big(\eta_4^{7}\wedge\eta_5^{7}+\eta_5^{7}\wedge\eta_6^{7}+\eta_6^{7}\wedge\eta_7^{7}\big)^2 \\
& = & \partial\eta_3^{7}\wedge\eta_3^{7}\wedge\big(2\, \eta_4^{7}\wedge\eta_5^{7}\wedge\eta_6^{7}\wedge\eta_7^{7}\big) \\
& = & 2\, \partial\eta_3^{7}\wedge\big(a_{33}\cdots a_{77}\, \varphi_3\wedge\varphi_4\wedge\varphi_5\wedge\varphi_6\wedge\varphi_7\big) \\
& = & 2\, \big(a_{33}\, \varphi_2\wedge\varphi_1\big)\wedge\big(a_{33}\cdots a_{77}\, \varphi_3\wedge\varphi_4\wedge\varphi_5\wedge\varphi_6\wedge\varphi_7\big)  \\
& = & -2\, a_{33} (\det A)\, 
\varphi_1\wedge\varphi_2\wedge \varphi_3\wedge\varphi_4\wedge\varphi_5\wedge\varphi_6\wedge\varphi_7\\
& \not= & 0.
\end{eqnarray*}

Suppose that the result holds for $l$, i.e. $\Gamma_l\wedge\partial \Gamma_l= c\,\varphi_{1}\wedge\cdots\wedge\varphi_{4l+3}$ with $c\not=0$. 
Recall that $\Gamma_l=\Phi_l^{4l+3}$ and $\Gamma_{l+1}=\Phi_{l+1}^{4(l+1)+3}$. We can write the last form as
$$
\Phi_{l+1}^{4(l+1)+3}= \Phi_{l}^{4(l+1)+3}\wedge \big(\eta_{4l+4}^{4l+7}\wedge\eta_{4l+5}^{4l+7}+\eta_{4l+5}^{4l+7}\wedge\eta_{4l+6}^{4l+7}+\eta_{4l+6}^{4l+7}\wedge\eta_{4l+7}^{4l+7}\big).
$$
Therefore,
\begin{eqnarray*}
\partial\Phi_{l+1}^{4(l+1)+3} & = & \partial\Phi_{l}^{4(l+1)+3}\wedge \big(\eta_{4l+4}^{4l+7}\wedge\eta_{4l+5}^{4l+7}+\eta_{4l+5}^{4l+7}\wedge\eta_{4l+6}^{4l+7}+\eta_{4l+6}^{4l+7}\wedge\eta_{4l+7}^{4l+7}\big) \\ 
&& - \Phi_{l}^{4(l+1)+3}\wedge \partial\big(\eta_{4l+4}^{4l+7}\wedge\eta_{4l+5}^{4l+7}+\eta_{4l+5}^{4l+7}\wedge\eta_{4l+6}^{4l+7}+\eta_{4l+6}^{4l+7}\wedge\eta_{4l+7}^{4l+7}\big),
\end{eqnarray*}
and using that $\Phi_{l}^{4(l+1)+3}\wedge\Phi_{l}^{4(l+1)+3}=0$, because $\Phi_{l}^{4(l+1)+3}$ is a form of odd degree, we get
\begin{eqnarray*}
\Gamma_{l+1}\wedge\partial \Gamma_{l+1} & = &  \Phi_{l+1}^{4(l+1)+3}\wedge \partial\Phi_{l+1}^{4(l+1)+3} \\
& = & 
\Phi_{l}^{4l+7}\wedge\partial\Phi_{l}^{4l+7}\wedge \big(\eta_{4l+4}^{4l+7}\wedge\eta_{4l+5}^{4l+7}+\eta_{4l+5}^{4l+7}\wedge\eta_{4l+6}^{4l+7}+\eta_{4l+6}^{4l+7}\wedge\eta_{4l+7}^{4l+7}\big)^2 \\
& = & \Phi_{l}^{4l+7}\wedge\partial\Phi_{l}^{4l+7}\wedge \big(2\, \eta_{4l+4}^{4l+7}\wedge\eta_{4l+5}^{4l+7}\wedge\eta_{4l+6}^{4l+7}\wedge\eta_{4l+7}^{4l+7}\big) \\
& = & 2\, a_{4l+4\, 4l+4}\cdots a_{4l+7\,4l+7}\, \Phi_{l}^{4l+7}\wedge\partial\Phi_{l}^{4l+7}\wedge \varphi_{4l+4}\wedge\varphi_{4l+5}\wedge\varphi_{4l+6}\wedge\varphi_{4l+7} \\
& = & 2\, a_{4l+4\, 4l+4}\cdots a_{4l+7\,4l+7}\, \Phi_{l}^{4l+3}\wedge\partial\Phi_{l}^{4l+3}\wedge \varphi_{4l+4}\wedge\varphi_{4l+5}\wedge\varphi_{4l+6}\wedge\varphi_{4l+7}, \end{eqnarray*}
where the last equality is due to the fact that any element in $\Phi_{l}^{4l+7}\wedge\partial\Phi_{l}^{4l+7}$ containing $\varphi_{4l+4}$, $\varphi_{4l+5}$, $\varphi_{4l+6}$ or $\varphi_{4l+7}$ is annihilated by the $(4,\,0)$-form $\varphi_{4l+4}\wedge\varphi_{4l+5}\wedge\varphi_{4l+6}\wedge\varphi_{4l+7}$.

Now, using that $\Phi_l^{4l+3}=\Gamma_l$ and the induction hypothesis, we conclude 
\begin{eqnarray*}
\Gamma_{l+1}\wedge\partial \Gamma_{l+1} & = & 2\, a_{4l+4\, 4l+4}\cdots a_{4l+7\,4l+7}\, \Gamma_{l}\wedge\partial \Gamma_{l}\wedge \varphi_{4l+4}\wedge\varphi_{4l+5}\wedge\varphi_{4l+6}\wedge\varphi_{4l+7}  \\ 
& = & 2\, c\, a_{4l+4\, 4l+4}\cdots a_{4l+7\,4l+7}\, \varphi_{1}\wedge\cdots\wedge\varphi_{4l+3}\wedge \varphi_{4l+4}\wedge\varphi_{4l+5}\wedge\varphi_{4l+6}\wedge\varphi_{4l+7} \\
& = & c'\, \varphi_{1}\wedge\cdots\wedge\varphi_{4l+7},
\end{eqnarray*}
for a non-zero constant $c'$. Thus, $\Gamma_{l+1}$ is a holomorphic $(p=2l+3)$-contact structure on the corresponding $(4l+7)$-dimensional complex manifold $X$. \hfill $\Box$

\begin{Obs}\label{Obs:p-contact-example} One can construct many examples in {\bf class I} by making different choices of invertible upper triangular matrices $A$, apart from the identity matrix. For instance, we can take $A$ defined by $a_{ij}=1$ for any $i\leq j$, which in complex dimension $7$ gives the following holomorphic $3$-contact structure 
\begin{eqnarray*}
\Gamma = \sum\limits_{3\leq i<j <k \leq 7}\varphi_{i}\wedge\varphi_{j}\wedge\varphi_{k},
\end{eqnarray*} 
in terms of the basis of holomorphic $(1,\,0)$-forms $\varphi_1,\dots , \varphi_7$.

\end{Obs}

\subsection{Examples induced by holomorphic $s$-symplectic manifolds}\label{subsection:examples_hol-s-symplectic} Recall the following classical notion: a {\it holomorphic symplectic} structure (or form) on a compact complex manifold $X$ of dimension $2s$ is a $d$-closed holomorphic $2$-form $\omega$ that is non-degenerate at every point. This means that $\omega\in C^\infty_{2,\,0}(X,\,\C)$ has the properties: \begin{eqnarray*}(i)\,\, d\omega = 0; \hspace{3ex} (ii)\,\,  \bar\partial\omega = 0; \hspace{3ex} (iii)\,\, \omega^s\neq 0 \hspace{2ex} \mbox{at every point of}\hspace{1ex} X.\end{eqnarray*} Since $\omega$ is of pure type, (i) implies (ii), but we wish to stress that, besides being holomorphic, $\omega$ is required to be $d$-closed.

\vspace{1ex}

We start by proposing the following generalisation of this classical notion.

\begin{Def}\label{Def:hol-s-symplectic} Let $X$ be a compact complex manifold with $\mbox{dim}_\C X = 2s$.

\vspace{1ex}
  
 $(1)$\, A {\bf holomorphic $s$-symplectic structure} on $X$ is a smooth $(s,\,0)$-form $\Omega\in C^\infty_{s,\,0}(X,\,\C)$ such that \begin{eqnarray*}(i)\,\,\bar\partial\Omega = 0; \hspace{3ex} \mbox{and}  \hspace{3ex} (ii)\,\, \Omega\wedge\Omega \neq 0 \hspace{2ex} \mbox{at every point of}\hspace{1ex} X.\end{eqnarray*} 

\vspace{1ex}
  
$(2)$\, We say that $X$ is a {\bf holomorphic $s$-symplectic manifold} if there exists a holomorphic $s$-symplectic structure on $X$.

\end{Def}  

Note that we do not require holomorphic $s$-symplectic structures $\Omega$ to be $d$-closed.

Moreover, property (ii) in their definition implies that $s$ is {\bf even} (since otherwise $\Omega\wedge\Omega$ would vanish identically). Thus, the complex dimension of any {\it holomorphic $s$-symplectic manifold} is a multiple of $4$. Furthermore, if $\Omega$ is an $s$-symplectic structure on $X$ with $\mbox{dim}_\C X = n = 2s$, the holomorphic $(n,\,0)$-form $u_\Omega:=\Omega\wedge\Omega$ is non-vanishing, so the canonical bundle $K_X$ is trivial. Thus, every holomorphic $s$-symplectic manifold $X$ is a Calabi-Yau manifold.  

If $s=2l$ with $l$ a positive integer, any {\it holomorphic symplectic} structure $\omega\in C^\infty_{2,\,0}(X,\,\C)$ on $X$ induces the {\it holomorphic $s$-symplectic structure} $\Omega:=\omega^s$.

\begin{The}\label{The:hol-s-symplectic_hol-p-contact} Let $G$ be a nilpotent complex Lie group with $\mbox{dim}_\C G = n = 2s = 4l$ and let $\Lambda$ be a co-compact lattice in $G$. Let $\{\varphi_1,\dots , \varphi_{4l}\}$ be a $\C$-basis of the dual vector space $\frg^\star$ of the Lie algebra $\frg$ of $G$. We still denote by $\{\varphi_1,\dots , \varphi_{4l}\}$ the induced $\C$-basis of $H^{1,\,0}_{\bar\partial}(Y,\,\C)$, where $Y=G/\Lambda$ is the induced quotient compact complex $n$-dimensional manifold.

\vspace{1ex}

$(1)$\, The $(s,\,0)$-form \begin{eqnarray*}\Omega:=\varphi_1\wedge\dots\wedge\varphi_{2l} + \varphi_{2l+1}\wedge\dots\wedge\varphi_{4l}\end{eqnarray*} is a {\bf holomorphic $s$-symplectic structure} on $Y$.

\vspace{1ex}

$(2)$\, Set $p:=2l+1$ and consider the $(2p+1)$-dimensional compact complex nilmanifold $X$ defined by a basis $\{\pi^\star\varphi_1,\dots , \pi^\star\varphi_{4l}, \varphi_{4l+1}, \varphi_{4l+2}, \varphi_{4l+3}\}$ of holomorphic $(1,\,0)$-forms whose first $4l$ members are the pullbacks under the natural projection $\pi:X\longrightarrow Y$ of the forms considered under $(1)$ and the three extra members satisfy the structure equations on $X$: \begin{eqnarray}\label{eqn:structure-eqn_hol-s-symplectic_hol-p-contact}\partial\varphi_{4l+1} = \partial\varphi_{4l+2} = 0 \hspace{3ex}  \mbox{and} \hspace{3ex} \partial\varphi_{4l+3} = \varphi_{4l+1}\wedge\varphi_{4l+2} + \pi^\star\sigma,\end{eqnarray} where $\sigma$ is any rational $d$-closed $(2,\,0)$-form on $Y$.

Then, the $(p,\,0)$-form \begin{eqnarray*}\Gamma:=\pi^\star\Omega\wedge\varphi_{4l+3}\end{eqnarray*} is a {\bf holomorphic $p$-contact structure} on $X$.

\end{The}  

\noindent {\it Proof.} $(1)$\, Since $\bar\partial\varphi_j = 0$ for every $j$, we get $\bar\partial\Omega = 0$. Moreover, \begin{eqnarray*}\Omega\wedge\Omega = 2\,\varphi_1\wedge\dots\wedge\varphi_{4l}\end{eqnarray*} and this is a non-vanishing $(n,\,0)$-form on $Y$. This proves the contention.

\vspace{1ex}

$(2)$\, Straightforward calculations yield: \begin{eqnarray*}\bar\partial\Gamma & = & \pi^\star\bar\partial\Omega\wedge\varphi_{4l+3} +  \pi^\star\Omega\wedge\bar\partial\varphi_{4l+3} = 0 \\
  \partial\Gamma & = & \pi^\star\partial\Omega\wedge\varphi_{4l+3} + \pi^\star\Omega\wedge\partial\varphi_{4l+3}\end{eqnarray*} hence  \begin{eqnarray*}\Gamma\wedge\partial\Gamma & = & \pi^\star\Omega\wedge\varphi_{4l+3}\wedge(\pi^\star\partial\Omega\wedge\varphi_{4l+3} + \pi^\star\Omega\wedge\partial\varphi_{4l+3}) \\
    & = & \varphi_{4l+3}\wedge(\varphi_{4l+1}\wedge\varphi_{4l+2} + \pi^\star\sigma)\wedge\pi^\star\Omega^2 \\
    & = & \varphi_{4l+1}\wedge\varphi_{4l+2}\wedge\varphi_{4l+3}\wedge\pi^\star\Omega^2 = 2\,\varphi_1\wedge\dots\wedge\varphi_{4l}\wedge\varphi_{4l+1}\wedge\varphi_{4l+2}\wedge\varphi_{4l+3}.\end{eqnarray*} This last form is a non-vanishing $(n,\,0)$-form on $X$, so the contention is now proven. Note that we used the equalities: $\bar\partial\varphi_{4l+3} = 0$, $\bar\partial\Omega = 0$ and, holding for bidegree reasons on $Y$, $\sigma\wedge\Omega^2 = 0$ (which implies $\pi^\star\sigma\wedge\pi^\star\Omega^2 = 0$ on $X$). 

Note that the existence of a co-compact lattice is guaranteed by $\sigma$ being rational. \hfill $\Box$

\vspace{3ex}

The form $\sigma$ was introduced in (\ref{eqn:structure-eqn_hol-s-symplectic_hol-p-contact}) in order to prevent the manifold $X$ from being the product $Y\times I^{(3)}$ of $Y$ with the Iwasawa manifold $I^{(3)}$ when $\sigma\neq 0$. In this way, we get further holomorphic $p$-contact manifolds $X$ from holomorphic $s$-symplectic manifolds $Y$ through the above construction besides the products $Y\times I^{(3)}$ that we get in the case $\sigma = 0$. The $4l+3$ structure equations satisfied by the holomorphic $(1,\,0)$-forms $\varphi_1,\dots , \varphi_{4l+3}$ (of which (\ref{eqn:structure-eqn_hol-s-symplectic_hol-p-contact}) are the last three ones) define a complex Lie group $H$ that admits a lattice when $\sigma$ is {\it rational}. When $\sigma = 0$, we have $H = G\times G_0$, where $G_0$ is the Heisenberg group of $3\times 3$-matrices whose quotient by a certain lattice is the Iwasawa manifold $I^{(3)}$. When $\sigma\neq 0$, $H$ need not be the product of $G$ by another Lie group.

\vspace{2ex}

In the following example, we exhibit holomorphic $s$-symplectic manifolds in all dimensions $4l$ with $l\geq 2$ that are not holomorphic symplectic.

\begin{Ex}\label{Ex:s-symplectic_non-hol-symplectic} Let $n = 2s = 4l$, with $l\geq 2$ a positive integer, and let $G$ be the nilpotent $n$-dimensional complex Lie group equipped with the complex structure defined by the following structure equations involving a basis of holomorphic $(1,\,0)$-forms $\varphi_1,\dots , \varphi_{4l}$:

  \vspace{1ex}

\hspace{5ex} $d\varphi_1 = \cdots =d\varphi_{4l-2} =d\varphi_{4l-1} = 0, \hspace{2ex} d\varphi_{4l} = \varphi_1\wedge\varphi_2 + \varphi_3\wedge\varphi_4+ \cdots + \varphi_{4l-3}\wedge\varphi_{4l-2}$.

  \vspace{1ex}

  \noindent The Lie group $G$ is the product of the complex $(4l-1)$-dimensional Heisenberg group $H_{\mathbb{C}}^{4l-1}$ by $\mathbb{C}$, so $G$ has a co-compact lattice $\Lambda$. On the compact quotient manifold $X$, by Nomizu's theorem, the de Rham class of any closed $(2,\,0)$-form $\omega$ on $X$ is represented by a left-invariant closed $(2,\,0)$-form $\tilde{\omega}$.

  Next, we will see that any such $\tilde{\omega}$ is degenerate, so $\{\omega\}_{dR}^{2l}=\{\tilde{\omega}^{2l}\}_{dR}=0$ and $X$ does not admit any holomorphic symplectic structure.

Any left-invariant $(2,\,0)$-form $\tilde{\omega}$ can be written as  
 $\tilde{\omega}=\sum\limits_{1\leq i<j\leq 4l} \lambda_{ij} \, \varphi_{i}\wedge\varphi_{j}$, for some $\lambda_{ij}\in \mathbb{C}$. From the structure equations we get

 \vspace{1ex}

\hspace{5ex} $d\tilde{\omega} = -\sum\limits_{1\leq i\leq 4l-1} \lambda_{i\, 4l} \, \varphi_{i}\wedge d\varphi_{4l} = -\bigg(\sum\limits_{1\leq i\leq 4l-2} \lambda_{i\, 4l} \, \varphi_{i}\bigg)\wedge d\varphi_{4l} 
- \lambda_{4l-1\, 4l} \, \varphi_{4l-1}\wedge d\varphi_{4l},$

 \vspace{1ex}
 
 \noindent which implies that $d\tilde{\omega} =0$ if and only if $\bigg(\sum\limits_{1\leq i\leq 4l-2} \lambda_{i\, 4l} \, \varphi_{i}\bigg)\wedge d\varphi_{4l}=0$ and $\lambda_{4l-1\, 4l} \, \varphi_{4l-1}\wedge d\varphi_{4l}=0$. The reason is that the differential $d\varphi_{4l}\in \bigwedge^2 W$, where $W$ is the $\C$-vector space $W=\langle \varphi_1,\ldots,\varphi_{4l-2} \rangle$.

   This immediately gives $\lambda_{4l-1\, 4l} =0$. Moreover, $d\varphi_{4l}=\varphi_1\wedge\varphi_2 + \cdots + \varphi_{4l-3}\wedge\varphi_{4l-2}$ is a non-degenerate 2-form on $W$, so the map $\bigwedge^1W\longrightarrow \bigwedge^3W$ is injective, which implies that $\sum\limits_{1\leq i\leq 4l-2} \lambda_{i\, 4l} \, \varphi_{i}=0$, i.e. $\lambda_{1\, 4l} =\cdots=\lambda_{4l-2\, 4l} =0$.

   Summing up, $\tilde{\omega}=\sum\limits_{1\leq i<j\leq 4l-1} \lambda_{ij} \, \varphi_{i}\wedge\varphi_{j}$, so $\tilde{\omega}^{2l}=0$.
 
\end{Ex}

\vspace{2ex}

In the following example we exhibit another construction of a holomorphic 3-contact manifold, which is slightly different from the one given in (2) of Theorem~\ref{The:hol-s-symplectic_hol-p-contact}. It can be thought of as the holomorphic analogue of the definition of a $G_2$-structure in seven real dimensions starting from three self-dual 2-forms on a 4-torus.

\begin{Ex}\label{Ex:G_2_analogue} Consider the complex 4-dimensional torus defined by four holomorphic $(1,\,0)$-forms $\varphi_1,\ldots,\varphi_4$ with 
$d\varphi_1=\cdots=d\varphi_4=0$. Let $\omega_1,\omega_2,\omega_3$ be the holomorphic symplectic forms given by

 \vspace{2ex}

\hspace{5ex} $\omega_1 = \varphi_1\wedge\varphi_2+\varphi_3\wedge\varphi_4,\quad
\omega_2 = \varphi_1\wedge\varphi_3-\varphi_2\wedge\varphi_4,\quad
\omega_3 = \varphi_1\wedge\varphi_4+\varphi_2\wedge\varphi_3.
$

 \vspace{2ex}
 
 \noindent These forms satisfy the conditions:

 \vspace{1ex}
 
 \noindent $\bullet$ $\omega_i^2=\omega_i\wedge\omega_i=2\, \varphi_1\wedge\varphi_2\wedge\varphi_3\wedge\varphi_4$, for $1\leq i\leq 3$;

 \vspace{1ex}
 
 \noindent $\bullet$ $\omega_i\wedge\omega_j=0$, for any $i\not= j$.

\vspace{1ex}

 Let us consider the $7$-dimensional compact complex nilmanifold $X$ defined by a basis $\{\varphi_1,\dots, \varphi_{7}\}$ of holomorphic $(1,\,0)$-forms whose first four members are the above closed forms and the three extra members satisfy the structure equations: 
\begin{eqnarray}\label{eqn:structure-g2t}
d\varphi_5=\omega_1, \hspace{4ex}
d\varphi_6=\omega_2, \hspace{4ex}
d\varphi_7=\omega_3.
\end{eqnarray} 

Define the $(3,\,0)$-form $\Gamma$ on $X$ by

 \vspace{2ex}

\hspace{20ex} $\Gamma= \varphi_5\wedge\omega_1+ \varphi_6\wedge\omega_2 +\varphi_7\wedge\omega_3 + \varphi_5\wedge\varphi_6\wedge\varphi_7$.

 \vspace{2ex}

Then $\delbar\Gamma=0$ and $\del\Gamma=\omega_1^2+\omega_2^2+\omega_3^2+ \omega_1\wedge\varphi_6\wedge\varphi_7 - \varphi_5\wedge\omega_2\wedge\varphi_7 + \varphi_5\wedge\varphi_6\wedge\omega_3$. 
Therefore,
\begin{eqnarray*}
\Gamma\wedge\partial\Gamma & = & \varphi_5\wedge\omega_1\wedge(\omega_1\wedge\varphi_6\wedge\varphi_7) + 
\varphi_6\wedge\omega_2\wedge(-\varphi_5\wedge\omega_2\wedge\varphi_7) +
\varphi_7\wedge\omega_3\wedge(\varphi_5\wedge\varphi_6\wedge\omega_3) \\ 
&& + 
\varphi_5\wedge\varphi_6\wedge\varphi_7\wedge(\omega_1^2+\omega_2^2+\omega_3^2) \\
    & = & 2 (\omega_1^2+\omega_2^2+\omega_3^2)\wedge\varphi_5\wedge\varphi_6\wedge\varphi_7 = 12\, \varphi_1\wedge\cdots\wedge\varphi_7.
 \end{eqnarray*} 
Consequently, $\Gamma\wedge\partial\Gamma\not=0$ at every point of $X$. Hence, $\Gamma$ is a holomorphic $3$-contact structure on $X$. 

Furthermore, $X$ does not admit any holomorphic contact structure.
Indeed, since $X$ is complex parallelisable, any holomorphic $1$-form $\gamma$ is left invariant, so $\gamma=\sum_{i=1}^{7} \lambda_i\,\varphi_i$ for some $\lambda_i\in\C$. Now, from \eqref{eqn:structure-g2t} we get $\partial \gamma\wedge \partial \gamma = 2 \big(\lambda_5^2+\lambda_6^2+\lambda_7^2\big)\varphi_1\wedge\varphi_2\wedge\varphi_3\wedge\varphi_4$, which implies $\gamma\wedge(\partial \gamma)^3=0$.
\end{Ex}

\vspace{2ex}

The following observation provides a direct way of defining higher-degree holomorphic contact manifolds starting from the examples given above throughout this section.

\begin{Prop}\label{Prop:products} Let $X$ be a holomorphic $p$-contact manifold and $Y$ a holomorphic $s$-symplectic manifold. Then, $Z=X\times Y$ is a holomorphic $(p+s)$-contact manifold.

\end{Prop}

\noindent {\it Proof.} Let $\Gamma_X$ be a holomorphic $p$-contact structure on $X$, i.e. $\delbar \Gamma_X=0$ and $\Gamma_X\wedge\del \Gamma_X\not=0$ at every point of $X$. Let $\Omega_Y$ be a holomorphic $s$-symplectic structure on $Y$, i.e. $\delbar \Omega_Y=0$ and $\Omega_Y\wedge \Omega_Y\not=0$ at every point of $Y$. Then, on the product complex manifold $Z$, the $(p+s)$-form $\Gamma=\Gamma_X\wedge\Omega_Y$ is holomorphic, and 
$\Gamma\wedge\del \Gamma=\Gamma_X\wedge\Omega_Y\wedge(\del \Gamma_X\wedge\Omega_Y - \Gamma_X\wedge\del\Omega_Y)=\Gamma_X\wedge\del \Gamma_X\wedge\Omega_Y\wedge\Omega_Y\not=0$ at every point of $Z$.
\hfill $\Box$

\section{The sheaves ${\cal F}_\Gamma$ and ${\cal G}_\Gamma$}\label{section:sheaves_F-G}

These two subsheaves of the holomorphic tangent sheaf ${\cal O}(T^{1,\,0}X)$ are naturally associated with the forms $\Gamma$ and $\partial\Gamma$ of a holomorphic $p$-contact structure.

\begin{Def}\label{Def:sheaves_F-G} Suppose $\Gamma\in C^\infty_{p,\,0}(X,\,\C)$ is a holomorphic $p$-contact structure on a compact complex manifold $X$ with $\mbox{dim}_\C X = n = 2p+1$.

\vspace{1ex}

(i)\, We let ${\cal F}_\Gamma$ be the sheaf of germs of holomorphic $(1,\,0)$-vector fields $\xi$ such that $\xi\lrcorner\Gamma=0$.

\vspace{1ex}

(ii)\, We let ${\cal G}_\Gamma$ be the sheaf of germs of holomorphic $(1,\,0)$-vector fields $\xi$ such that $\xi\lrcorner\partial\Gamma=0$.

\end{Def}

Considering the holomorphic vector bundle morphism \begin{eqnarray*}T^{1,\,0}X\longrightarrow\Lambda^{p-1,\,0}T^\star X, \hspace{3ex} \xi\longmapsto\xi\lrcorner\Gamma,\end{eqnarray*} the ${\cal O}_X$-module ${\cal F}_\Gamma$ is the kernel of the induced morphism between the locally free sheaves associated with these vector bundles. In particular, ${\cal F}_\Gamma$ is a {\it coherent} analytic subsheaf of ${\cal O}(T^{1,\,0}X)$ (since the kernel of a morphism of coherent sheaves is coherent -- see e.g. [Dem97, II, Theorem 3.13]). Moreover, ${\cal F}_\Gamma$ is {\it torsion-free} since any coherent subsheaf of a torsion-free sheaf is torsion-free. Meanwhile, it is standard that every torsion-free coherent sheaf is locally free outside an analytic subset of codimension $\geq 2$. (See e.g. [Kob87, Corollary 5.5.15].) Thus, there exists an analytic subset $\Sigma_\Gamma\subset X$ with $\mbox{codim}\,\Sigma_\Gamma\geq 2$ such that ${\cal F}_\Gamma$ is locally free on $X\setminus\Sigma_\Gamma$.

Similarly, by considering the holomorphic vector bundle morphism \begin{eqnarray*}T^{1,\,0}X\longrightarrow\Lambda^{p,\,0}T^\star X, \hspace{3ex} \xi\longmapsto\xi\lrcorner\partial\Gamma,\end{eqnarray*} and its sheaf counterpart, we deduce that the latter's kernel, ${\cal G}_\Gamma$, is a {\it torsion-free coherent} analytic subsheaf of ${\cal O}(T^{1,\,0}X)$. In particular, ${\cal G}_\Gamma$ is locally free on $X\setminus\Xi_{\partial\Gamma}$, where $\Xi_{\partial\Gamma}\subset X$ is an analytic subset such that $\mbox{codim}\,\Xi_{\partial\Gamma}\geq 2$.

\vspace{1ex}

Since ${\cal F}_\Gamma$ and ${\cal G}_\Gamma$ are kernels of holomorphic vector bundle morphisms, they are holomorphic subbundles of $T^{1,\,0}X$ (i.e. the analytic subsets $\Sigma_\Gamma$, resp. $\Xi_{\partial\Gamma}$, are empty) if and only if the maps $T^{1,\,0}_xX\ni\xi\longmapsto\xi\lrcorner\Gamma\in\Lambda^{p-1,\,0}T^\star_xX$, respectively $T^{1,\,0}_xX\ni\xi\longmapsto\xi\lrcorner\partial\Gamma\in\Lambda^{p,\,0}T^\star_xX$, are of ranks independent of the point $x\in X$. 

\begin{Prop}\label{Prop:F-G-sheaves_properties} Let $X$ be a compact complex manifold with $\mbox{dim}_\C X = n = 2p+1$. Suppose there exists a {\bf holomorphic $p$-contact} structure $\Gamma\in C^\infty_{p,\,0}(X,\,\C)$ on $X$.

\vspace{1ex}

(i)\, The sum ${\cal F}_\Gamma\oplus{\cal G}_\Gamma\subset{\cal O}(T^{1,\,0}X)$ is direct.

\vspace{1ex}

(ii)\, If there exists a {\bf holomorphic contact} structure $\eta\in C^\infty_{1,\,0}(X,\,\C)$, the sheaf of germs ${\cal F}_\eta$ of holomorphic $(1,\,0)$-vector fields $\xi$ such that $\xi\lrcorner\eta=0$ and the sheaf of germs ${\cal G}_\eta$ of holomorphic $(1,\,0)$-vector fields $\xi$ such that $\xi\lrcorner\partial\eta=0$ are {\bf locally free} of respective ranks $n-1$ and $1$, while $T^{1,\,0}X$ has a direct-sum splitting \begin{eqnarray}\label{eqn:contact_direct-sum-splitting}T^{1,\,0}X = {\cal F}_\eta\oplus{\cal G}_\eta\end{eqnarray} with ${\cal F}_\eta$ and ${\cal G}_\eta$ viewed as holomorphic vector subbundles of $T^{1,\,0}X$.

  Furthermore, the holomorphic line bundle ${\cal G}_\eta$ is {\bf trivial} and $\eta$ is a non-vanishing global holomorphic section of its dual ${\cal G}^\star_\eta$.

 \vspace{1ex}

 (iii) Suppose that ${\cal F}_\Gamma$ and ${\cal G}_\Gamma$ are {\bf locally free} and induce a direct-sum splitting $T^{1,\,0}X = {\cal F}_\Gamma\oplus{\cal G}_\Gamma$ into vector subbundles. For every $s=0,\dots , n$ and every $\varphi\in C^\infty_{0,\,s}(X,\,T^{1,\,0}X)$, let $\varphi = \varphi_\Gamma + \varphi_{\partial\Gamma}$ be the induced splitting of $\varphi$ into pieces $\varphi_\Gamma, \varphi_{\partial\Gamma}\in C^\infty_{0,\,s}(X,\,T^{1,\,0}X)$ such that $\varphi_\Gamma\lrcorner\Gamma = 0$ and $\varphi_{\partial\Gamma}\lrcorner\partial\Gamma = 0$.

 Fix an arbitrary $q\in\{1,\dots , n\}$.

\vspace{1ex}

(a)\, For any $\theta\in C^\infty_{0,\,q}(X,\,T^{1,\,0}X)$ such that $\theta\in\mbox{Im}\,\bar\partial$ and $\theta\lrcorner\Gamma = 0$ and for any $\varphi\in C^\infty_{0,\,q-1}(X,\,T^{1,\,0}X)$ such that $\bar\partial\varphi = \theta$, one has $\bar\partial\varphi_\Gamma = \theta$ and $\bar\partial\varphi_{\partial\Gamma} = 0$.

 \vspace{1ex}

 (b)\, For any $\theta\in C^\infty_{0,\,q}(X,\,T^{1,\,0}X)$ such that $\theta\in\mbox{Im}\,\bar\partial$ and $\theta\lrcorner\partial\Gamma = 0$ and for any $\varphi\in C^\infty_{0,\,q-1}(X,\,T^{1,\,0}X)$ such that $\bar\partial\varphi = \theta$, one has $\bar\partial\varphi_{\partial\Gamma} = \theta$ and $\bar\partial\varphi_\Gamma = 0$.

 \vspace{1ex}

 (c)\, For any $\theta\in C^\infty_{0,\,q}(X,\,T^{1,\,0}X)$, the following equivalences hold: \begin{eqnarray}\label{eqn:equivalences_closed_pieces-F-G}\theta\in\ker\bar\partial & \iff & \theta_\Gamma\in\ker\bar\partial \hspace{2ex}\mbox{and}\hspace{2ex} \theta_{\partial\Gamma}\in\ker\bar\partial \\
   \label{eqn:equivalences_exact_pieces-F-G}\theta\in\mbox{Im}\,\bar\partial & \iff & \theta_\Gamma\in\mbox{Im}\,\bar\partial \hspace{2ex}\mbox{and}\hspace{2ex} \theta_{\partial\Gamma}\in\mbox{Im}\,\bar\partial.\end{eqnarray}

\vspace{1ex}

(d)\, For every $s=0,\dots , n$, let $C^\infty_{0,\,s}(X,\,{\cal F}_\Gamma):=\bigg\{\varphi\in C^\infty_{0,\,s}(X,\,T^{1,\,0}X)\,\mid\,\varphi\lrcorner\Gamma = 0\bigg\}$ and  $C^\infty_{0,\,s}(X,\,{\cal G}_\Gamma):=\bigg\{\varphi\in C^\infty_{0,\,s}(X,\,T^{1,\,0}X)\,\mid\,\varphi\lrcorner\partial\Gamma = 0\bigg\}$.

The $\C$-vector spaces: \begin{eqnarray}\label{eqn:H_Gamma-del-Gamma_def}
  \nonumber H^{0,\,q}_\Gamma(X,\,T^{1,\,0}X) & := & \frac{\ker\bigg(\bar\partial:C^\infty_{0,\,q}(X,\,{\cal F}_\Gamma)\longrightarrow C^\infty_{0,\,q+1}(X,\,{\cal F}_\Gamma)\bigg)}{\mbox{Im}\,\bigg(\bar\partial:C^\infty_{0,\,q-1}(X,\,{\cal F}_\Gamma)\longrightarrow C^\infty_{0,\,q}(X,\,{\cal F}_\Gamma)\bigg)}, \\
  H^{0,\,q}_{\partial\Gamma}(X,\,T^{1,\,0}X) & := & \frac{\ker\bigg(\bar\partial:C^\infty_{0,\,q}(X,\,{\cal G}_\Gamma)\longrightarrow C^\infty_{0,\,q+1}(X,\,{\cal G}_\Gamma)\bigg)}{\mbox{Im}\,\bigg(\bar\partial:C^\infty_{0,\,q-1}(X,\,{\cal G}_\Gamma)\longrightarrow C^\infty_{0,\,q}(X,\,{\cal G}_\Gamma)\bigg)}\end{eqnarray} are well defined and the $\C$-linear map \begin{eqnarray}\label{eqn:H_bar_Gamma-del-Gamma_splitting}\nonumber \Phi: H^{0,\,q}_{\bar\partial}(X,\,T^{1,\,0}X) & \longrightarrow & H^{0,\,q}_\Gamma(X,\,T^{1,\,0}X)\oplus H^{0,\,q}_{\partial\Gamma}(X,\,T^{1,\,0}X) \\
 \{\theta\}_{\bar\partial} & \longmapsto & \bigg([\theta_\Gamma]_{\bar\partial},\,[\theta_{\partial\Gamma}]_{\bar\partial}\bigg)\end{eqnarray} is well defined and an {\bf isomorphism}.

\end{Prop}  

\noindent {\it Proof.} (i)\, Suppose $\xi\in({\cal F}_\Gamma)_x\cap({\cal G}_\Gamma)_x$ for some $x\in X$. Then, $\xi\lrcorner\Gamma = 0$ and $\xi\lrcorner\partial\Gamma = 0$ on some neighbourhood $U$ of $x$. This implies that $\xi\lrcorner u_\Gamma = (\xi\lrcorner\Gamma)\wedge\partial\Gamma + (-1)^p\,\Gamma\wedge(\xi\lrcorner\partial\Gamma) = 0$, which in turn implies that $\xi = 0$ on $U$ since the contraction of the Calabi-Yau form $u_\Gamma$ is an isomorphism.

\vspace{1ex}

(ii)\, The ${\cal O}_X$-module ${\cal F}_\eta$ is the kernel of the morphism of locally free sheaves induced by the holomorphic vector bundle morphism \begin{eqnarray*}T^{1,\,0}X\longrightarrow\Lambda^{0,\,0}T^\star X, \hspace{3ex} \xi\longmapsto\xi\lrcorner\eta.\end{eqnarray*} Since the rank of $\Lambda^{0,\,0}T^\star X$ is $1$ and this map is non-identically zero at every point $x\in X$ (due to the hypothesis $\eta\wedge(\partial\eta)^p\neq 0$ everywhere), this map must be surjective at every $x\in X$. Thus, its kernel has dimension $n-1$ at every point, which implies that the sheaf ${\cal F}_\eta$ is locally free of rank $n-1$.

Meanwhile, the ${\cal O}_X$-module ${\cal G}_\eta$ is the kernel of the morphism of locally free sheaves induced by the holomorphic vector bundle morphism \begin{eqnarray*}T^{1,\,0}X\longrightarrow\Lambda^{1,\,0}T^\star X, \hspace{3ex} \xi\longmapsto\xi\lrcorner\partial\eta.\end{eqnarray*}

Now, by the argument that gave (i) we infer that the sum $({\cal F}_\eta)_x\oplus({\cal G}_\eta)_x\subset{\cal O}(T^{1,\,0}X)_x$ is direct, hence the dimension of $({\cal G}_\eta)_x$ is either $0$ or $1$, for every $x\in X$. Meanwhile, if there were a point $x\in X$ such that the dimension of $({\cal G}_\eta)_x$ equals $0$, the contraction map $T^{1,\,0}_xX\ni\xi\longmapsto\xi\lrcorner\partial\eta\in\Lambda^{1,\,0}T^\star_xX$ would be bijective. This would amount to the alternating form $(\partial\eta)_x:T^{1,\,0}_xX\times T^{1,\,0}_xX\longrightarrow\C$ being non-degenerate, hence to its (necessarily skew-symmetric) matrix $A$ (e.g. with respect to the basis of $T^{1,\,0}_xX$ induced by a given system of local coordinates $z_1,\dots , z_n$ centred at $x$) being invertible. From $^tA = -A$ we would then get $\det A = (-1)^n\,\det A$, which would imply, thanks to $\det A\neq 0$, that $n$ is even, a contradiction.  

We conclude that the dimension of $({\cal G}_\eta)_x$ equals $1$ at every point $x\in X$, hence the sheaf ${\cal G}_\eta$ is locally free of rank $1$ and (\ref{eqn:contact_direct-sum-splitting}) holds.

The definition of ${\cal F}_\eta$ means that its fibre at every point $x\in X$ is precisely the kernel of the $\C$-linear map $\eta_x:T^{1,\,0}_xX = ({\cal F}_\eta)_x\oplus({\cal G}_\eta)_x\longrightarrow\C$. Thus, for every $x\in X$, the restriction of $\eta_x$ to $({\cal G}_\eta)_x$ is a $\C$-linear isomorphism onto $\C$. This amounts to $\eta$ being a non-vanishing global holomorphic section of ${\cal G}^\star_\eta$, which implies the triviality of this holomorphic line bundle. Then, its dual ${\cal G}_\eta$ is trivial as well. 

\vspace{1ex}

(iii)\, Let $\theta\in C^\infty_{0,\,q}(X,\,T^{1,\,0}X)$ such that $\theta\in\mbox{Im}\,\bar\partial$ and let $\varphi\in C^\infty_{0,\,q-1}(X,\,T^{1,\,0}X)$ such that $\bar\partial\varphi = \theta$. This implies, thanks also to $\bar\partial\Gamma = 0$ and to $\bar\partial(\partial\Gamma) = 0$, that \begin{eqnarray}\label{xi_Gamma_del-Gamma_del-bar}\nonumber\bar\partial(\varphi_\Gamma\lrcorner\partial\Gamma) & = & \bar\partial(\varphi\lrcorner\partial\Gamma) = (\bar\partial\varphi)\lrcorner\partial\Gamma = \theta\lrcorner\partial\Gamma, \\
 \bar\partial(\varphi_{\partial\Gamma}\lrcorner\Gamma) & = & \bar\partial(\varphi\lrcorner\Gamma) = (\bar\partial\varphi)\lrcorner\Gamma = \theta\lrcorner\Gamma.\end{eqnarray}

\vspace{1ex}

(a)\, We get: \begin{eqnarray*}(\bar\partial\varphi_\Gamma)\lrcorner u_\Gamma = \bar\partial(\varphi_\Gamma\lrcorner u_\Gamma) = (-1)^{pq}\,\bar\partial\bigg(\Gamma\wedge(\varphi_\Gamma\lrcorner\partial\Gamma)\bigg) = (-1)^{p(q+1)}\,\Gamma\wedge\bar\partial(\varphi_\Gamma\lrcorner\partial\Gamma) = (-1)^{p(q+1)}\,\Gamma\wedge(\theta\lrcorner\partial\Gamma),\end{eqnarray*} where the last equality follows from the former line in (\ref{xi_Gamma_del-Gamma_del-bar}).

Suppose now that $\theta\lrcorner\Gamma = 0$. Then, $\theta\lrcorner u_\Gamma = (-1)^{p(q+1)}\,\Gamma\wedge(\theta\lrcorner\partial\Gamma)$. When combined with the above equality, this implies $(\bar\partial\varphi_\Gamma)\lrcorner u_\Gamma = \theta\lrcorner u_\Gamma$, which amounts to $\bar\partial\varphi_\Gamma = \theta$ (cf. the Calabi-Yau isomorphism (\ref{eqn:C_Y-isomorphism_scalar})), as claimed. Since $\bar\partial\varphi = \theta$, this implies that $\bar\partial\varphi_{\partial\Gamma} = 0$.

\vspace{1ex}

(b)\, When $\theta\lrcorner\partial\Gamma = 0$, $\theta\lrcorner u_\Gamma = (\theta\lrcorner\Gamma)\wedge\partial\Gamma$. This yields the last equality below: \begin{eqnarray*}(\bar\partial\varphi_{\partial\Gamma})\lrcorner u_\Gamma = \bar\partial(\varphi_{\partial\Gamma}\lrcorner u_\Gamma) = \bar\partial\bigg((\varphi_{\partial\Gamma}\lrcorner\Gamma)\wedge\partial\Gamma\bigg) = \bar\partial(\varphi_{\partial\Gamma}\lrcorner\Gamma)\wedge\partial\Gamma = (\theta\lrcorner\Gamma)\wedge\partial\Gamma = \theta\lrcorner u_\Gamma,\end{eqnarray*} where the last but one equality follows from the latter line in (\ref{xi_Gamma_del-Gamma_del-bar}). Thanks to the Calabi-Yau isomorphism (\ref{eqn:C_Y-isomorphism_scalar}), this amounts to $\bar\partial\varphi_{\partial\Gamma} = \theta$ which further implies $\bar\partial\varphi_\Gamma = 0$, as claimed.

\vspace{1ex}

(c)\, $\bullet$ {\it Proof of equivalence (\ref{eqn:equivalences_closed_pieces-F-G}).} If $\bar\partial\theta_\Gamma = \bar\partial\theta_{\partial\Gamma} = 0$, then $\bar\partial\theta = \bar\partial\theta_\Gamma + \bar\partial\theta_{\partial\Gamma} = 0$.

Conversely, suppose that $\bar\partial\theta = 0$. Then \begin{eqnarray*}0 = (\bar\partial\theta)\lrcorner\Gamma = \bar\partial(\theta\lrcorner\Gamma) = \bar\partial(\theta_{\partial\Gamma}\lrcorner\Gamma) = (\bar\partial\theta_{\partial\Gamma})\lrcorner\Gamma.\end{eqnarray*} Meanwhile, $(\bar\partial\theta_{\partial\Gamma})\lrcorner\partial\Gamma = \bar\partial(\theta_{\partial\Gamma}\lrcorner\partial\Gamma) = 0$ because $\theta_{\partial\Gamma}\lrcorner\partial\Gamma = 0$. We conclude that $(\bar\partial\theta_{\partial\Gamma})\lrcorner u_\Gamma = 0$, hence $\bar\partial\theta_{\partial\Gamma} = 0$. Since we also have $\bar\partial\theta = 0$, we also get $\bar\partial\theta_\Gamma = 0$.

\vspace{1ex}

$\bullet$ {\it Proof of equivalence (\ref{eqn:equivalences_exact_pieces-F-G}).}
If $\theta_\Gamma,\, \theta_{\partial\Gamma}\in\mbox{Im}\,\bar\partial$, then $\theta = \theta_\Gamma + \theta_{\partial\Gamma}\in\mbox{Im}\,\bar\partial$.

Conversely, suppose there exists $\varphi\in C^\infty_{0,\,q-1}(X,\,T^{1,\,0}X)$ such that $\theta = \bar\partial\varphi$. This is equivalent to \begin{eqnarray*}A_\Gamma:=\theta_\Gamma - \bar\partial\varphi_\Gamma = \bar\partial\varphi_{\partial\Gamma} - \theta_{\partial\Gamma}.\end{eqnarray*}

Now, note that $A_\Gamma\lrcorner\Gamma = -\bar\partial(\varphi_\Gamma\lrcorner\Gamma) = 0$ since $\theta_\Gamma\lrcorner\Gamma = 0$ and $\varphi_\Gamma\lrcorner\Gamma = 0$, while $A_\Gamma\lrcorner\partial\Gamma = \bar\partial(\varphi_{\partial\Gamma}\lrcorner\partial\Gamma) = 0$ since $\theta_{\partial\Gamma}\lrcorner\partial\Gamma = 0$ and $\varphi_{\partial\Gamma}\lrcorner\partial\Gamma = 0$. We conclude that $A_\Gamma\lrcorner u_\Gamma = 0$, hence $A_\Gamma = 0$. This means that \begin{eqnarray*}\theta_\Gamma = \bar\partial\varphi_\Gamma \hspace{3ex}\mbox{and}\hspace{3ex} \theta_{\partial\Gamma} = \bar\partial\varphi_{\partial\Gamma}.\end{eqnarray*}

\vspace{1ex}

(d)\, To show the well-definedness of the map $\bar\partial:C^\infty_{0,\,q}(X,\,{\cal F}_\Gamma)\longrightarrow C^\infty_{0,\,q+1}(X,\,{\cal F}_\Gamma)$, we have to check that whenever $\varphi\in C^\infty_{0,\,q}(X,\,T^{1,\,0}X)$ has the property $\varphi\lrcorner\Gamma = 0$, it also satisfies $(\bar\partial\varphi)\lrcorner\Gamma = 0$. This follows at once from $0 = \bar\partial(\varphi\lrcorner\Gamma) = (\bar\partial\varphi)\lrcorner\Gamma$.

Meanwhile, for every class $[\theta]_{\bar\partial}\in H^{0,\,q}_{\bar\partial}(X,\,T^{1,\,0}X)$ and every representative $\theta$ thereof, one has $\bar\partial\theta = 0$, hence (c) ensures that $\bar\partial\theta_\Gamma = 0$ and  $\bar\partial\theta_{\partial\Gamma} = 0$. Since $\theta_\Gamma\in C^\infty_{0,\,q}(X,\,{\cal F}_\Gamma)$ and $\theta_{\partial\Gamma}\in C^\infty_{0,\,q}(X,\,{\cal G}_\Gamma)$, we get well-defined classes $[\theta_\Gamma]_{\bar\partial}\in H^{0,\,q}_\Gamma(X,\,T^{1,\,0}X)$ and $[\theta_{\partial\Gamma}]_{\bar\partial}\in H^{0,\,q}_{\partial\Gamma}(X,\,T^{1,\,0}X)$.

Moreover, $\Phi$ is independent of the choice of representative $\theta$ of the
class $[\theta]_{\bar\partial}\in H^{0,\,q}_{\bar\partial}(X,\,T^{1,\,0}X)$. Indeed, whenever $\theta\in\mbox{Im}\,\bar\partial$, (c) and (a) ensure that $[\theta_\Gamma]_{\bar\partial} = 0$ in $H^{0,\,q}_\Gamma(X,\,T^{1,\,0}X)$ and $[\theta_{\partial\Gamma}]_{\bar\partial} = 0$ in $H^{0,\,q}_{\partial\Gamma}(X,\,T^{1,\,0}X)$. This finishes the proof of the well-definedness of $\Phi$.

The injectivity and the surjectivity of $\Phi$ follow at once from the definitions.     \hfill $\Box$

\subsection{Further properties of the sheaves ${\cal F}_\Gamma$ and ${\cal G}_\Gamma$}\label{subsection:further-prop_F-G}

\vspace{3ex}

The following remark means that ${\cal G}_\Gamma$ defines a {\it (possibly singular) holomorphic foliation}. As seen above, its (possibly empty) set of singularities is $\Xi_{\partial\Gamma}\subset X$. 

\begin{Obs}\label{Obs:G_Gamma_integrability} The subsheaf ${\cal G}_\Gamma$ of ${\cal O}(T^{1,\,0}X)$ is {\bf integrable} in the sense that $[{\cal G}_\Gamma,\,{\cal G}_\Gamma]\subset{\cal G}_\Gamma$, where $[\,\cdot\,,\,\cdot\,]$ is the Lie bracket of $T^{1,\,0}X$.

\end{Obs}

\noindent {\it Proof.} Since $d\partial\Gamma = 0 = d(d\Gamma)$ (because $\bar\partial\Gamma = 0$), the Cartan formula reads: \begin{eqnarray}\label{eqn:Cartan-formula_G_del-Gamma}\nonumber 0= d(\partial\Gamma)(\xi_0,\dots ,\xi_{p+1}) & = & \sum\limits_{j=0}^{p+1}(-1)^j\,\xi_j\cdot(\partial\Gamma)(\xi_0,\dots ,\widehat\xi_j,\dots , \xi_{p+1}) \\
  & +  & \sum\limits_{0\leq j<k\leq p+1}(-1)^{j+k}\,(\partial\Gamma)([\xi_j,\,\xi_k],\,\xi_0,\dots,\widehat\xi_j,\dots ,\widehat\xi_k,\dots , \xi_{p+1})\end{eqnarray} for all $(1,\,0)$-vector fields $\xi_0,\dots , \xi_{p+1}$.

If two among $\xi_0,\dots ,\xi_{p+1}$, say $\xi_l$ and $\xi_s$ for some $l<s$, are (local) sections of ${\cal G}_\Gamma$, then:

\vspace{1ex}

$\bullet$ all the terms in the first sum on the r.h.s. of (\ref{eqn:Cartan-formula_G_del-Gamma}) vanish;

 \vspace{1ex}

 $\bullet$ all the terms in the second sum on the r.h.s. of (\ref{eqn:Cartan-formula_G_del-Gamma}) vanish, except possibly the term $$(-1)^{l+s}\,(\partial\Gamma)([\xi_l,\,\xi_s],\,\xi_0,\dots,\widehat\xi_l,\dots ,\widehat\xi_s,\dots , \xi_{p+1}).$$

 \vspace{1ex}

 Thus, this last term must vanish as well for all $(1,\,0)$-vector fields $\xi_0,\dots,\widehat\xi_l,\dots ,\widehat\xi_s,\dots , \xi_{p+1}$. This means that $[\xi_l,\,\xi_s]\lrcorner\partial\Gamma = 0$, which amounts to $[\xi_l,\,\xi_s]$ being a (local) section of ${\cal G}_\Gamma$, for all pairs $\xi_l,\,\xi_s$ of (local) sections of ${\cal G}_\Gamma$.

 This proves the integrability of ${\cal G}_\Gamma$.  \hfill $\Box$

\vspace{2ex}

Based on the considerations in $\S$\ref{subsection:Lie}, we now introduce directional analogues, suited to our setting, of the page-$1$-$\partial\bar\partial$- and $\partial\bar\partial$-properties.

\begin{Def}\label{Def:dd-bar_F-G_directional} Let $(X,\,\Gamma)$ be a compact complex {\bf holomorphic $p$-contact manifold} with $\mbox{dim}_\C X = n =2p+1$.

\vspace{1ex}

$(1)$\, We say that $(X,\,\Gamma)$ is a {\bf partially ${\cal F}_\Gamma$-directional page-$1$-$\partial\bar\partial$-manifold} if for every $\theta\in C^\infty_{0,\,1}(X,\,{\cal F}_\Gamma)$ such that $\bar\partial\theta = 0$ and $\theta\lrcorner\partial\Gamma\in{\cal Z}_2^{p,\,1}(X)$, there exists $\xi\in C^\infty(X,\,{\cal F}_\Gamma)$ such that $\partial(\theta\lrcorner\partial\Gamma) = \partial\bar\partial(\xi\lrcorner\partial\Gamma)$.

\vspace{1ex}

$(2)$\, We say that $(X,\,\Gamma)$ is a {\bf partially ${\cal G}_\Gamma$-directional page-$1$-$\partial\bar\partial$-manifold} if for every $\theta\in C^\infty_{0,\,1}(X,\,{\cal G}_\Gamma)$ such that $\bar\partial\theta = 0$ and $\theta\lrcorner\Gamma\in{\cal Z}_2^{p-1,\,1}(X)$, there exists $\xi\in C^\infty(X,\,{\cal G}_\Gamma)$ such that $\partial(\theta\lrcorner\Gamma) = \partial\bar\partial(\xi\lrcorner\Gamma)$.

\vspace{1ex}

$(3)$\, We say that $(X,\,\Gamma)$ is a {\bf partially vertically $\partial\bar\partial$-manifold} if:

\vspace{1ex}

(a)\, for every $\theta\in C^\infty_{0,\,1}(X,\,{\cal G}_\Gamma)$ such that $\partial(\theta\lrcorner\Gamma)\in\ker\bar\partial$, there exists $\xi\in C^\infty(X,\,{\cal G}_\Gamma)$ such that $\partial(\theta\lrcorner\Gamma) = \partial\bar\partial(\xi\lrcorner\Gamma)$;

\vspace{1ex}

(b)\,  for all constantly vertical $\theta_1,\dots\theta_{2k}\in C^\infty_{0,\,1}(X,\,T^{1,\,0}X)$ and $a_1,\dots , a_k\in\C$, the following implication holds: \begin{eqnarray*}\partial\bigg(\sum\limits_{s=1}^k a_s\,\theta_{2s-1}\lrcorner(\theta_{2s}\lrcorner\Gamma)\bigg)\in\ker\bar\partial \implies \partial\bigg(\sum\limits_{s=1}^k a_s\,\theta_{2s-1}\lrcorner(\theta_{2s}\lrcorner\Gamma)\bigg)\in\mbox{Im}\,\bar\partial.\end{eqnarray*}

\end{Def}

Note that properties $(1)$ and $(2)$ in the above definition are directional versions of property (\ref{eqn:ddbar_im_Z-2}), itself one of the aspects of the {\it page-$1$-$\partial\bar\partial$-property} introduced and studied in [PSU20a] and [PSU20b], hence the adverb ``partially''. Meanwhile, property $(3)$ is a version in the ${\cal G}_\Gamma$-directions of the inclusions $\mbox{Im}\,\partial\cap\ker\bar\partial\subset\mbox{Im}\,(\partial\bar\partial)$ and $\mbox{Im}\,\partial\cap\ker\bar\partial\subset\mbox{Im}\,\bar\partial$, themselves two of the inclusions that constitute the classical {\it $\partial\bar\partial$-property}. In all three cases, the $\partial\bar\partial$-potential whose existence is asserted therein has the same properties as the $\partial$-potential given beforehand. 

\vspace{2ex}

We now refine the $\C$-vector spaces $H^{0,\,1}_\Gamma(X,\,T^{1,\,0}X)$ and $H^{0,\,1}_{\partial\Gamma}(X,\,T^{1,\,0}X)$ introduced in (\ref{eqn:H_Gamma-del-Gamma_def}) by defining subspaces that can be seen as directional analogues of the $E_2$-spaces on the second page of the Fr\"olicher spectral sequence: \begin{eqnarray}\label{eqn:E_2_directional}\nonumber E^{0,\,1}_2(X,\,{\cal F}_\Gamma) & := & \bigg\{[\theta]_{\bar\partial}\in H^{0,\,1}_\Gamma(X,\,T^{1,\,0}X)\,\mid\,\bar\partial\theta = 0,\hspace{1ex} \theta\lrcorner\Gamma = 0, \hspace{1ex} \theta\lrcorner\partial\Gamma\in{\cal Z}_2^{p,\,1}(X)\bigg\} \\
  E^{0,\,1}_2(X,\,{\cal G}_\Gamma) & := & \bigg\{[\theta]_{\bar\partial}\in H^{0,\,1}_{\partial\Gamma}(X,\,T^{1,\,0}X)\,\mid\,\bar\partial\theta = 0,\hspace{1ex} \theta\lrcorner\partial\Gamma = 0, \hspace{1ex} \theta\lrcorner\Gamma\in{\cal Z}_2^{p-1,\,1}(X)\bigg\}.\end{eqnarray}

Finally, we show that, under the appropriate assumption of the type spelt out in Definition \ref{Def:dd-bar_F-G_directional}, this directional analogue of the $E_2$-cohomology consists only of classes that are representable by {\it constantly horizontal}, respectively {\it constantly vertical}, vector-valued forms.

\begin{Prop}\label{Prop:E_2_directional} Let $(X,\,\Gamma)$ be a compact complex {\bf holomorphic $p$-contact manifold} with $\mbox{dim}_\C X = n = 2p+1$.

\vspace{1ex}

$(1)$\, Suppose that $(X,\,\Gamma)$ is a {\bf partially ${\cal F}_\Gamma$-directional page-$1$-$\partial\bar\partial$-manifold}.

Then, for every $\theta\in C^\infty_{0,\,1}(X,\,{\cal F}_\Gamma)$ such that $\bar\partial\theta = 0$ and $\theta\lrcorner\partial\Gamma\in{\cal Z}_2^{p,\,1}(X)$, there exists $\psi\in C^\infty_{0,\,1}(X,\,{\cal F}_\Gamma)$ such that $\bar\partial\psi = 0$, $[\psi]_{\bar\partial} = [\theta]_{\bar\partial}\in H^{0,\,1}_\Gamma(X,\,T^{1,\,0}X)$ and $\partial(\psi\lrcorner\partial\Gamma) = 0$.

In other words, every class $[\theta]_{\bar\partial}\in E^{0,\,1}_2(X,\,{\cal F}_\Gamma)$ has a {\bf constantly horizontal} representative $\psi$.

\vspace{1ex}

$(2)$\, Suppose that $(X,\,\Gamma)$ is a {\bf partially ${\cal G}_\Gamma$-directional page-$1$-$\partial\bar\partial$-manifold}.

Then, for every $\theta\in C^\infty_{0,\,1}(X,\,{\cal G}_\Gamma)$ such that $\bar\partial\theta = 0$ and $\theta\lrcorner\Gamma\in{\cal Z}_2^{p-1,\,1}(X)$, there exists $\psi\in C^\infty_{0,\,1}(X,\,{\cal G}_\Gamma)$ such that $\bar\partial\psi = 0$, $[\psi]_{\bar\partial} = [\theta]_{\bar\partial}\in H^{0,\,1}_{\partial\Gamma}(X,\,T^{1,\,0}X)$ and $\partial(\psi\lrcorner\Gamma) = 0$.

In other words, every class $[\theta]_{\bar\partial}\in E^{0,\,1}_2(X,\,{\cal G}_\Gamma)$ has a {\bf constantly vertical} representative $\psi$.

\end{Prop}  

\noindent {\it Proof.} $(1)$\, The partial ${\cal F}_\Gamma$-directional page-$1$-$\partial\bar\partial$-hypothesis ensures the existence of $\xi\in C^\infty(X,\,{\cal F}_\Gamma)$ such that $$\partial\bigg((\theta - \bar\partial\xi)\lrcorner\partial\Gamma\bigg) = 0.$$ Setting $\psi := \theta - \bar\partial\xi$ proves the contention.

\vspace{1ex}

$(2)$\, The partial ${\cal G}_\Gamma$-directional page-$1$-$\partial\bar\partial$-hypothesis ensures the existence of $\xi\in C^\infty(X,\,{\cal G}_\Gamma)$ such that $$\partial\bigg((\theta - \bar\partial\xi)\lrcorner\Gamma\bigg) = 0.$$ Setting $\psi := \theta - \bar\partial\xi$ proves the contention.  \hfill $\Box$

\vspace{2ex}

 As for ${\cal F}_\Gamma$, it need not be integrable, but it may have an analogous property in cohomology (cf. (i) in the next definition) which may coexist with a maximal non-integrability property at the pointwise level (cf. (ii) below) that is more in tune with the idea behind a $p$-contact structure.

 \begin{Def}\label{Def:F_Gamma_coh-integrability} Let $X$ be a compact complex manifold with $\mbox{dim}_\C X = n = 2p+1$. Suppose there exists a holomorphic $p$-contact structure $\Gamma\in C^\infty_{p,\,0}(X,\,\C)$ on $X$.

\vspace{1ex}

(i)\, We say that the sheaf ${\cal F}_\Gamma$ is {\bf cohomologically integrable in bidegree $(0,\,1)$} if for all pairs of $\bar\partial$-cohomology classes $[\varphi]_{\bar\partial}, [\psi]_{\bar\partial}\in H^{0,\,1}_{\bar\partial}(X,\,T^{1,\,0}X)$ such that $[\varphi\lrcorner\Gamma]_{\bar\partial} = [\psi\lrcorner\Gamma]_{\bar\partial} = 0\in H^{p-1,\,1}_{\bar\partial}(X,\,\C)$ we have:  \begin{eqnarray*}\bigg[[\varphi,\,\psi]\lrcorner\Gamma\bigg]_{\bar\partial}= 0\in H^{p-1,\,2}_{\bar\partial}(X,\,\C),
\end{eqnarray*} 
where $[\varphi,\,\psi]$ is the standard Lie bracket recalled in (\ref{eqn:bracket_0pq-vector_def}) of the forms $\varphi,\psi\in C^\infty_{0,\,1}(X,\,T^{1,\,0}X)$.

\vspace{1ex}

(ii)\, We say that the sheaf ${\cal F}_\Gamma$ is {\bf constantly maximally non-integrable in bidegree $(0,\,1)$} if for all pairs of {\bf constantly horizontal} forms $\varphi, \psi\in C^\infty_{0,\,1}(X,\,T^{1,\,0}X)$, the form $[\varphi,\,\psi]\in C^\infty_{0,\,2}(X,\,T^{1,\,0}X)$ 
obtained by Lie bracket 
is {\bf vertical} (i.e. $[\varphi,\,\psi]\in C^\infty_{0,\,2}(X,\,{\cal G}_\Gamma)$, or equivalently $[\varphi,\,\psi]\lrcorner\partial\Gamma = 0$). 

 \end{Def}

 \section{A structure theorem}\label{section:structure-theorem} In Theorem \ref{The:hol-s-symplectic_hol-p-contact}, we constructed holomorphic $p$-contact manifolds $X$ from holomorphic $s$-symplectic manifolds $Y$ such that $X$ is fibred over $Y$ and both are nilmanifolds. In this section, we give a generalised converse of Theorem \ref{The:hol-s-symplectic_hol-p-contact}: we will show that, under certain assumptions, any holomorphic $p$-contact manifold $X$ arises as a fibration over a holomorphic $s$-symplectic manifold $Y$. For reasons explained in the Introduction, we assume the existence of a fibration whose fibres are the leaves of a given holomorphic foliation.

\begin{The}\label{Thm:structure_Iwasawa-s-symplectic} Let $\pi:X\longrightarrow Y$ be a surjective holomorphic submersion between compact complex manifolds with $\mbox{dim}_\C X = n = 2p+1 = 4l +3$ and $\mbox{dim}_\C Y = 4l = 2s$. Suppose that:

\vspace{1ex}

(a)\, $X$ has a {\bf holomorphic $p$-contact structure} $\Gamma\in C^\infty_{p,\,0}(X,\,\C)$;

\vspace{1ex}

(b)\, $T^{1,\,0}X = {\cal E}\oplus{\cal H}$, where ${\cal E}$ and ${\cal H}$ are holomorphic subbundles of $T^{1,\,0}X$ such that ${\cal H}$ is {\bf Frobenius integrable} and the leaves of the foliation it induces on $X$ are the fibres of $\pi:X\longrightarrow Y$;

\vspace{1ex}

(c)\, there exist holomorphic vector fields $\eta_1, \eta_2, \eta_3\in H^0(X,\,{\cal H})$ that globally trivialise ${\cal H}$ such that $\eta_1, \eta_2\in H^0(X,\,{\cal F}_\Gamma)$, $\psi_3\wedge\Gamma = 0$ (where $\psi_1,\,\psi_2,\,\psi_3\in H^{1,\,0}_{\bar\partial}(X,\,\C)$ is the global frame of ${\cal H}^\star$ dual to the global frame $\{\eta_1,\,\eta_2,\,\eta_3\}$ of ${\cal H}$) and the following relations are satisfied: \begin{eqnarray}\label{eqn:eta-s_relations}[\eta_1,\,\eta_2] = -[\eta_2,\,\eta_1] = \eta_3 \hspace{3ex}\mbox{and} \hspace{3ex} [\eta_j,\,\eta_k] = 0 \hspace{2ex} \mbox{for all}\hspace{1ex} (j,\,k)\notin\{(1,\,2),\,(2,\,1)\};\end{eqnarray} 
  
\vspace{1ex}

Then, there exists a {\bf holomorphic $s$-symplectic structure} $\Omega\in C^\infty_{p-1,\,0}(Y,\,\C) = C^\infty_{2l,\,0}(Y,\,\C) = C^\infty_{s,\,0}(Y,\,\C)$ on $Y$ such that \begin{eqnarray*}\Gamma = \pi^\star\Omega\wedge\psi_3 \hspace{3ex}\mbox{on}\hspace{1ex} X.\end{eqnarray*}

\end{The}

\noindent {\it Proof.} It proceeds in a number of steps.

\hspace{1ex}

$\bullet$ {\it Step $1$.} We define a form $\widetilde\Omega$ on $X$ that will subsequently be shown to be the pullback under $\pi$ of the holomorphic $s$-symplectic form $\Omega$ on $Y$ whose existence is claimed in the statement.

\begin{Lem}\label{Lem:Omega-tilde_def} The following equality of $(p-1,\,0)$-forms on $X$ holds: \begin{eqnarray*}\eta_1\lrcorner(\eta_2\lrcorner\partial\Gamma) = \eta_3\lrcorner\Gamma.\end{eqnarray*}

\end{Lem}

\noindent {\it Proof of Lemma.} The definition of $L^{1,\,0}_{\eta_2}$ yields: $\partial(\eta_2\lrcorner\Gamma) = L^{1,\,0}_{\eta_2}(\Gamma) - \eta_2\lrcorner\partial\Gamma$. Contracting by $\eta_1$, this gives: \begin{eqnarray*}\eta_1\lrcorner(\eta_2\lrcorner\partial\Gamma) = \eta_1\lrcorner L^{1,\,0}_{\eta_2}(\Gamma) - \eta_1\lrcorner\partial(\eta_2\lrcorner\Gamma).\end{eqnarray*}

\noindent Meanwhile, the definition of $L^{1,\,0}_{\eta_1}$ yields: \begin{eqnarray*}\partial\bigg(\eta_1\lrcorner(\eta_2\lrcorner\Gamma)\bigg) = L^{1,\,0}_{\eta_1}(\eta_2\lrcorner\Gamma) - \eta_1\lrcorner\partial(\eta_2\lrcorner\Gamma).\end{eqnarray*} Equating the value of $\eta_1\lrcorner\partial(\eta_2\lrcorner\Gamma)$ across the last two equalities, we get: \begin{eqnarray*}\eta_1\lrcorner(\eta_2\lrcorner\partial\Gamma) = \eta_1\lrcorner L^{1,\,0}_{\eta_2}(\Gamma) - L^{1,\,0}_{\eta_1}(\eta_2\lrcorner\Gamma) + \partial\bigg(\eta_1\lrcorner(\eta_2\lrcorner\Gamma)\bigg).\end{eqnarray*}

On the other hand, the general commutation property (a) of (\ref{eqn:L_commutations}) yields the first equality below: \begin{eqnarray*}[\eta_1,\,\eta_2]\lrcorner\Gamma = [\eta_1\lrcorner\cdot,\,L^{1,\,0}_{\eta_2}](\Gamma) = \eta_1\lrcorner L^{1,\,0}_{\eta_2}(\Gamma) - L^{1,\,0}_{\eta_2}(\eta_1\lrcorner\Gamma).\end{eqnarray*}

\noindent Equating the value of $\eta_1\lrcorner L^{1,\,0}_{\eta_2}(\Gamma)$ across the last two equalities, we get: \begin{eqnarray}\label{Lem:Omega-tilde_def_proof}\eta_1\lrcorner(\eta_2\lrcorner\partial\Gamma) = [\eta_1,\,\eta_2]\lrcorner\Gamma + L^{1,\,0}_{\eta_2}(\eta_1\lrcorner\Gamma) -  L^{1,\,0}_{\eta_1}(\eta_2\lrcorner\Gamma) + \partial\bigg(\eta_1\lrcorner(\eta_2\lrcorner\Gamma)\bigg).\end{eqnarray}

Now, $[\eta_1,\,\eta_2] = \eta_3$ (by assumption), while $\eta_1\lrcorner\Gamma = \eta_2\lrcorner\Gamma = 0$ (since $\eta_1$ and $\eta_2$ are supposed to be sections of ${\cal F}_\Gamma$), so (\ref{Lem:Omega-tilde_def_proof}) amounts to the contention. \hfill $\Box$

\vspace{2ex}

We now give a name to the form that is either of the two equal forms involved in Lemma \ref{Lem:Omega-tilde_def}:  \begin{eqnarray}\label{Def:Omega-tilde_s-symplectic}\widetilde\Omega: = \eta_1\lrcorner(\eta_2\lrcorner\partial\Gamma) = \eta_3\lrcorner\Gamma\in C^\infty_{p-1,\,0}(X,\,\C) = C^\infty_{2l,\,0}(X,\,\C) = C^\infty_{s,\,0}(X,\,\C).\end{eqnarray}

Since, by hypothesis, $\eta_1$, $\eta_2$ and $\eta_3$ are {\it holomorphic} $(1,\,0)$-vector fields and $\Gamma$ is a {\it holomorphic} $(p,\,0)$-form, we get \begin{eqnarray}\label{eqn::Omega-tilde_holomorphic}\bar\partial\widetilde\Omega = 0.\end{eqnarray}

\hspace{1ex}

$\bullet$ {\it Step $2$.} We now observe some of the properties of the form $\widetilde\Omega$ defined in (\ref{Def:Omega-tilde_s-symplectic}).

\begin{Lem}\label{Lem:Omega-tilde_properties} The following equalities hold on $X$: \begin{eqnarray*}\eta_1\lrcorner(\eta_2\lrcorner\partial\widetilde\Omega) = \eta_3\lrcorner\widetilde\Omega = 0.\end{eqnarray*}

\end{Lem}

\noindent {\it Proof.} Since $\widetilde\Omega = \eta_3\lrcorner\Gamma$, applying $\partial$ and using the definition of $L^{1,\,0}_{\eta_3}$, we get: \begin{eqnarray*}\partial\widetilde\Omega = L^{1,\,0}_{\eta_3}(\Gamma) -  \eta_3\lrcorner\partial\Gamma.\end{eqnarray*} Hence, taking the double contraction by $\eta_1$ and $\eta_2$, we get: \begin{eqnarray}\label{Lem:Omega-tilde_properties_proof_1}\nonumber\eta_1\lrcorner(\eta_2\lrcorner\partial\widetilde\Omega) & = & \eta_1\lrcorner\bigg(\eta_2\lrcorner L^{1,\,0}_{\eta_3}(\Gamma)\bigg) - \eta_1\lrcorner\bigg(\eta_2\lrcorner(\eta_3\lrcorner\partial\Gamma)\bigg) \\
  & = & \eta_1\lrcorner\bigg(L^{1,\,0}_{\eta_3}(\eta_2\lrcorner\Gamma) + [\eta_2,\,\eta_3]\lrcorner\Gamma\bigg) - \eta_1\lrcorner\bigg(\eta_2\lrcorner(\eta_3\lrcorner\partial\Gamma)\bigg) = - \eta_1\lrcorner\bigg(\eta_2\lrcorner(\eta_3\lrcorner\partial\Gamma)\bigg),\end{eqnarray} where the second equality followed from the general commutation property (a) of (\ref{eqn:L_commutations}), while the third equality followed from $\eta_2\lrcorner\Gamma = 0$ and $[\eta_2,\,\eta_3] = 0$ (both of which are part of the hypotheses).

In order to understand the last term above, we start by taking $\partial$ in the equality $\eta_2\lrcorner\Gamma = 0$ to get $0 = \partial(\eta_2\lrcorner\Gamma) = L^{1,\,0}_{\eta_2}(\Gamma) - \eta_2\lrcorner\partial\Gamma$. Thus, $\eta_2\lrcorner\partial\Gamma = L^{1,\,0}_{\eta_2}(\Gamma)$, so contracting by $\eta_1$ we get: \begin{eqnarray*}\eta_1\lrcorner(\eta_2\lrcorner\partial\Gamma) = \eta_1\lrcorner L^{1,\,0}_{\eta_2}(\Gamma) = L^{1,\,0}_{\eta_2}(\eta_1\lrcorner\Gamma) + [\eta_1,\,\eta_2]\lrcorner\Gamma = \eta_3\lrcorner\Gamma,\end{eqnarray*} where the second equality followed from the general commutation property (a) of (\ref{eqn:L_commutations}), while the last equality followed from $\eta_1\lrcorner\Gamma = 0$ and $[\eta_1,\,\eta_2] = \eta_3$.

Thus, further contracting by $\eta_3$ the above equality $\eta_1\lrcorner(\eta_2\lrcorner\partial\Gamma) = \eta_3\lrcorner\Gamma$ and using the vanishing of any form after a double contraction by a same vector field ($\eta_3$ in this case), we get \begin{eqnarray}\label{Lem:Omega-tilde_properties_proof_2}\eta_3\lrcorner\bigg(\eta_1\lrcorner(\eta_2\lrcorner\partial\Gamma)\bigg) = 0.\end{eqnarray}

Putting together (\ref{Lem:Omega-tilde_properties_proof_1}) and (\ref{Lem:Omega-tilde_properties_proof_2}), we infer the vanishing of $\eta_1\lrcorner(\eta_2\lrcorner\partial\widetilde\Omega)$. Meanwhile, we know from (\ref{Def:Omega-tilde_s-symplectic}) that $\widetilde\Omega = \eta_1\lrcorner(\eta_2\lrcorner\partial\Gamma)$, so (\ref{Lem:Omega-tilde_properties_proof_2}) expresses the vanishing of $\eta_3\lrcorner\widetilde\Omega$. The proof is complete. \hfill $\Box$

\begin{Lem}\label{Lem:Omega-tilde_product} Recall that $u_\Gamma = \Gamma\wedge\partial\Gamma$ is a nowhere vanishing holomorphic $(n,\,0)$-form on $X$.

  The following equality holds on $X$: \begin{eqnarray*}\widetilde\Omega\wedge\widetilde\Omega = \eta_1\lrcorner\bigg(\eta_2\lrcorner(\eta_3\lrcorner u_\Gamma)\bigg).\end{eqnarray*} In particular, $\widetilde\Omega\wedge\widetilde\Omega$ is a nowhere vanishing holomorphic $(n-3,\,0)$-form on $X$.

\end{Lem}

\noindent {\it Proof.} Using (\ref{Def:Omega-tilde_s-symplectic}), we get: \begin{eqnarray*}\widetilde\Omega\wedge\widetilde\Omega = \widetilde\Omega\wedge(\eta_3\lrcorner\Gamma) = \eta_3\lrcorner(\widetilde\Omega\wedge\Gamma) = \eta_3\lrcorner\bigg(\bigg(\eta_1\lrcorner(\eta_2\lrcorner\partial\Gamma)\bigg)\wedge\Gamma\bigg) = \eta_3\lrcorner\bigg(\eta_1\lrcorner(\eta_2\lrcorner(\partial\Gamma\wedge\Gamma)\bigg)\bigg),\end{eqnarray*} where we have used the identities $\eta_3\lrcorner\widetilde\Omega = 0$ (see Lemma \ref{Lem:Omega-tilde_properties}), $\widetilde\Omega = \eta_1\lrcorner(\eta_2\lrcorner\partial\Gamma)$ (see (\ref{Def:Omega-tilde_s-symplectic})) and $\eta_1\lrcorner\Gamma = \eta_2\lrcorner\Gamma = 0$ (see the hypotheses of Theorem \ref{Thm:structure_Iwasawa-s-symplectic}). This proves the former contention.

The latter contention follows from the former and from the hypotheses of Theorem \ref{Thm:structure_Iwasawa-s-symplectic}.    \hfill $\Box$

\hspace{1ex}

$\bullet$ {\it Step $3$.} We now define a holomorphic $s$-symplectic form $\Omega$ on $Y$ such that $\pi^\star\Omega = \widetilde\Omega$. The expression of $\widetilde\Omega$ involving $\partial\Gamma$ in (\ref{Def:Omega-tilde_s-symplectic}) will be used.

Thanks to our assumptions, for every $y\in Y$ and every $x\in X_y:=\pi^{-1}(y)\subset X$, the differential map of $\pi$ at $x$ \begin{eqnarray*}d_x\pi : T^{1,\,0}_xX\longrightarrow T^{1,\,0}_{\pi(x)}Y = T^{1,\,0}_yY\end{eqnarray*} is surjective and its kernel is $T^{1,\,0}_xX_y = {\cal H}_x$. In particular, its restriction \begin{eqnarray*}(d_x\pi)_{|{\cal E}_x} : {\cal E}_x\longrightarrow T^{1,\,0}_{\pi(x)}Y = T^{1,\,0}_yY\end{eqnarray*} is an isomorphism of $\C$-vector spaces.

    Thus, for every $y\in Y$, every $\nu_y\in T^{1,\,0}_yY$ and every $x\in X_y$, there exists a unique $\xi_x\in{\cal E}_x\subset T^{1,\,0}_xX$ such that $(d_x\pi)(\xi_x) = \nu_y$. Since $\pi$ is a holomorphic submersion, the vector field \begin{eqnarray*}X_y\ni x\longmapsto\xi_x = (d_x\pi)^{-1}(\nu_y)\in{\cal E}_x\subset T^{1,\,0}_xX\end{eqnarray*} obtained in this way is holomorphic.

      In other words, for every $y\in Y$ and every $\nu_y\in T^{1,\,0}_yY$, there exists a unique holomorphic vector field $\xi\in H^0(X_y,\, {\cal E}_{|X_y})\subset H^0(X_y,\, (T^{1,\,0}X)_{|X_y})$ of type $(1,\,0)$ on the compact complex manifold $X_y$ such that \begin{eqnarray*}(d_x\pi)(\xi_x) = \nu_y, \hspace{5ex} \mbox{for all}\hspace{1ex} x\in X_y.\end{eqnarray*}

      Using this information, we define an $(s,\,0)$-form $\Omega$ (recall that $s=p-1$) on $Y$ by \begin{eqnarray}\label{eqn:Omega_s-symplectic_def_Y}\Omega_y(\nu_{1,\,y},\dots , \nu_{s,\,y}):=\widetilde\Omega_x(\xi_{1,\,x},\dots , \xi_{s,\,x}),\end{eqnarray} for all $y\in Y$, $\nu_{1,\,y},\dots , \nu_{s,\,y}\in T^{1,\,0}_yY$ and $x\in X_y$, where $\xi_{1,\,x},\dots , \xi_{s,\,x}\in{\cal E}_x\subset T^{1,\,0}_xX$ are the unique vectors such that $(d_x\pi)(\xi_{j,\,x}) = \nu_{j,\,y}$ for all $j\in\{1,\dots , s\}$.

      We now show that definition (\ref{eqn:Omega_s-symplectic_def_Y}) is correct. This means that, for every $y\in Y$, the complex number $\widetilde\Omega_x(\xi_{1,\,x},\dots , \xi_{s,\,x})$ is independent of the choice of $x\in X_y$. Indeed, the $C^\infty$ function $\xi_s\lrcorner\dots\lrcorner\xi_2\lrcorner\xi_1\lrcorner\widetilde\Omega$ on $X_y$ has the property: \begin{eqnarray*}\bar\partial\bigg(\xi_s\lrcorner\dots\lrcorner\xi_2\lrcorner\xi_1\lrcorner\widetilde\Omega\bigg) = -\xi_s\lrcorner\bar\partial\bigg(\xi_{s-1}\lrcorner\dots\lrcorner\xi_2\lrcorner\xi_1\lrcorner\widetilde\Omega\bigg) = \dots = (-1)^s\,\xi_s\lrcorner\xi_{s-1}\lrcorner\dots\lrcorner\xi_1\lrcorner\bar\partial\widetilde\Omega =0\end{eqnarray*} since $\bar\partial\xi_j = 0$ for all $j\in\{1,\dots , s\}$ (see above) and $\bar\partial\widetilde\Omega =0$ (see (\ref{eqn::Omega-tilde_holomorphic})). Consequently, this function is holomorphic on the compact complex manifold $X_y$, hence constant.

We now easily infer the following

\begin{Lem}\label{Lem:Omega_s-symplectic_properties} The form $\Omega\in C^\infty_{s,\,0}(Y,\,\C)$ defined in (\ref{eqn:Omega_s-symplectic_def_Y}) has the properties:  \begin{eqnarray}\label{eqn:Omega_s-symplectic_properties}(a)\hspace{1ex}\pi^\star\Omega = \widetilde\Omega; \hspace{3ex} (b)\hspace{1ex} \bar\partial\Omega = 0 \hspace{3ex}\mbox{and}\hspace{3ex} (c)\hspace{1ex}\Omega\wedge\Omega\neq 0 \hspace{2ex}\mbox{everywhere on}\hspace{1ex} Y.\end{eqnarray}

  In particular, $\Omega$ is a {\bf holomorphic $s$-symplectic form} on $Y$.

\end{Lem}
  
\noindent {\it Proof of Lemma.} Conclusions (a) and (b) follow right away from the definition (\ref{eqn:Omega_s-symplectic_def_Y}) of $\Omega$ and from $\bar\partial\widetilde\Omega = 0$ (cf. (\ref{eqn::Omega-tilde_holomorphic})), while (c) follows from (a) and from Lemma \ref{Lem:Omega-tilde_product}. \hfill $\Box$

\hspace{1ex}

$\bullet$ {\it Step $4$.} It remains to show that $\Gamma$ splits as product of $\widetilde\Omega = \pi^\star\Omega$ by $\psi_3$.

Since $\widetilde\Omega = \eta_3\lrcorner\Gamma$ (cf. (\ref{Def:Omega-tilde_s-symplectic})), we get: \begin{eqnarray*}\widetilde\Omega\wedge\psi_3 = (\eta_3\lrcorner\Gamma)\wedge\psi_3 = \eta_3\lrcorner(\Gamma\wedge\psi_3) + \Gamma = \Gamma,\end{eqnarray*} where the last two equalities followed from $\eta_3\lrcorner\psi_3 = 1$, respectively $\Gamma\wedge\psi_3 = 0$.

The proof of Theorem \ref{Thm:structure_Iwasawa-s-symplectic} is complete.  \hfill $\Box$

\vspace{2ex}

  To further put our results in their context, let us say that in Foreman's setting of [For00], where $(X,\,\eta)$ is a compact holomorphic contact manifold of dimension $n=2p+1$ and $\pi:X\longrightarrow Y$ is a holomorphic fibration whose fibres are the leaves of ${\cal G}_\eta$, the form $\partial\eta\in C^\infty_{2,\,0}(X,\,\C)$ can be ``copied'' as a form $\Omega$ on $Y$ as we did for the form $\widetilde\Omega$ (defined on $X$) that was ``copied'' as the form $\Omega$ (defined on $Y$) in Step $3$ of the proof of Theorem \ref{Thm:structure_Iwasawa-s-symplectic} (see (\ref{eqn:Omega_s-symplectic_def_Y})). Then, the form $\Omega\in C^\infty_{2,\,0}(Y,\,\C)$ can be shown to be holomorphic symplectic on $Y$ and to have the property $\pi^\star\Omega = d\eta$ on $X$ as we did in Lemma \ref{Lem:Omega_s-symplectic_properties} in our case.

 \section{Essential horizontal small deformations}\label{section:contact-def} In this section, we define a type of small deformations of the complex structure of $X$ that preserve a given holomorphic $p$-contact structure $\Gamma$. We then prove a certain form of unobstructedness.

\vspace{1ex}

\subsection{Meaning of unobstructedness for arbitrary small deformations}\label{subsection:meaning-unobstructedness} We start by recalling a few general, well-known facts in order to fix the notation and the terminology.

Let $X$ be an $n$-dimensional compact complex manifold. A cohomology class $[\theta_1]_{\bar\partial}\in H^{0,\,1}_{\bar\partial}(X,\,T^{1,\,0}X)$ is said to induce a family $(J_t)_{t\in D}$ of small deformations of the complex structure $J_0$ of $X_0:=X$ in the direction of $[\theta_1]_{\bar\partial}$ (where $D$ is a small disc about $0$ in $\C$) if there exists a representative $\psi_1\in C^\infty_{0,\,1}(X,\,T^{1,\,0}X)$ of $[\theta_1]_{\bar\partial}$ and a sequence $(\psi_j)_{j\geq 2}$ of vector-valued forms $\psi_j\in C^\infty_{0,\,1}(X,\,T^{1,\,0}X)$ such that, for all $t\in\C$ sufficiently close to $0$, the power series \begin{eqnarray*}\psi(t)=\psi_1\,t + \psi_2\,t^2 + \dots + \psi_\nu\,t^\nu + \dots\end{eqnarray*} converges (in some topology) and defines a vector-valued form $\psi(t)\in C^\infty_{0,\,1}(X,\,T^{1,\,0}X)$ that satisfies the {\it integrability condition} \begin{eqnarray}\label{eqn:integrability-condition}\bar\partial\psi(t) = \frac{1}{2}\,[\psi(t), \,\psi(t)], \hspace{6ex} t\in D.\end{eqnarray}

  This integrability condition is easily seen to be equivalent to $\bar\partial\psi_1=0$ (always satisfied since $\psi_1$ represents a $\bar\partial$-cohomology class) holding simultaneously with the following sequence of equations: \begin{equation}\label{eqn:integrability-condition_nu}\bar\partial\psi_\nu = \frac{1}{2}\,\sum\limits_{\mu=1}^{\nu-1}[\psi_\mu,\,\psi_{\nu-\mu}] \hspace{3ex} (\mbox{Eq.}\,\,(\nu)), \hspace{3ex} \nu\geq 2. \end{equation}

  The integrability condition (\ref{eqn:integrability-condition}) means that the almost complex structure $J_t$ defined on $X$ by $\psi(t)$ via the operator \begin{eqnarray*}\bar\partial_t\simeq\bar\partial_0 - \psi(t)\end{eqnarray*} is {\it integrable} in the sense that $\bar\partial_t^2 = 0$. (Here, $\bar\partial_0$ is the Cauchy-Riemann operator associated with $J_0$ and the meaning of $\simeq$ is that a $C^\infty$ function $f$ on an open subset of $X$ is holomorphic w.r.t. $J_t$ in the sense that $\bar\partial_t f = 0$ if and only if $(\bar\partial_0 - \psi(t))\,f = 0$.) This amounts to $J_t$ being a complex structure on $X$. In this case, one regards $J_t$ as a deformation of $J_0$.

  Geometrically, the original cohomology class \begin{eqnarray*}[\theta_1]_{\bar\partial} = \bigg[\frac{\partial\psi(t)}{\partial t}_{|t=0}\bigg]_{\bar\partial}\in H^{0,\,1}(X,\,T^{1,\,0}X)\end{eqnarray*} arises, if the integrability condition (\ref{eqn:integrability-condition}) is satisfied, as the tangent vector at $t=0$ to the complex curve $(J_t)_{t\in D}$. Therefore, one regards $[\theta_1]_{\bar\partial}$ as a deformation to order $1$ of $J_0$.

  It is standard to say that the {\it Kuranishi family of $X$ is unobstructed} (to arbitrary order) if every cohomology class $[\theta_1]_{\bar\partial}\in H^{0,\,1}_{\bar\partial}(X,\,T^{1,\,0}X)$ induces a family $(J_t)_{t\in D}$ of small deformations of the complex structure $J_0$ of $X_0:=X$ in the direction of $[\theta_1]_{\bar\partial}$. This is equivalent to the solvability of all the equations in (\ref{eqn:integrability-condition_nu}) for every $[\theta_1]_{\bar\partial}\in H^{0,\,1}_{\bar\partial}(X,\,T^{1,\,0}X)$ and for some representative $\psi_1\in C^\infty_{0,\,1}(X,\,T^{1,\,0}X)$ of the class $[\theta_1]_{\bar\partial}$.

  \subsection{Unobstructedness for essential horizontal deformations}\label{subsection:unobstructedness_arbitrary-order} The point of view that we take on deformations in this paper originates in the following

\begin{Def}\label{Def:essential-horizontal-def} Let $X$ be a compact complex manifold with $\mbox{dim}_\C X = n = 2p+1$. Suppose there exists a holomorphic $p$-contact structure $\Gamma\in C^\infty_{p,\,0}(X,\,\C)$ on $X$.

  (i)\, The $\C$-vector subspace \begin{eqnarray*} E^{0,\,1}_2(X,\,{\cal F}_\Gamma)\subset H^{0,\,1}_\Gamma(X,\,T^{1,\,0}X)\end{eqnarray*} defined in (\ref{eqn:E_2_directional}) is called the space of (infinitesimal) {\bf essential horizontal deformations} of $X$.

  \vspace{1ex}

(ii)\, We say that the essential horizontal deformations of $X$ are {\bf unobstructed} if every class $[\theta_1]_{\bar\partial}\in E^{0,\,1}_2(X,\,{\cal F}_\Gamma)$ has a representative $\psi_1\in C^\infty_{0,\,1}(X,\,T^{1,\,0}X)$ such that all the equations $(Eq.\, (\nu))$ in (\ref{eqn:integrability-condition_nu}) admit solutions $\psi_\nu\in C^\infty_{0,\,1}(X,\,T^{1,\,0}X)$ with the property $\psi_\nu\lrcorner u_\Gamma\in\ker\partial$ for every $\nu\geq 2$.

\end{Def}

Note that there is a well-defined linear map from the space $E^{0,\,1}_2(X,\,{\cal F}_\Gamma)\subset H^{0,\,1}_\Gamma(X,\,T^{1,\,0}X)$ of {\it essential horizontal deformations} of $X$ to the space $H^{0,\,1}_{[\Gamma],\,def}(X,\,T^{1,\,0}X)\subset H^{0,\,1}_{\bar\partial}(X,\,T^{1,\,0}X)$ of {\it $p$-contact deformations} of $X$ introduced in [KPU25, $\S8.2$]. This map sends every class $[\theta]_{\bar\partial}$ in the former space to the class of the same representative in the latter space. However, these two classes are different since the quotients defining the spaces $H^{0,\,1}_\Gamma(X,\,T^{1,\,0}X)$ and $H^{0,\,1}_{\bar\partial}(X,\,T^{1,\,0}X)$ are computed differently. Moreover, thanks to (iii)(a) of Proposition \ref{Prop:F-G-sheaves_properties}, said map $E^{0,\,1}_2(X,\,{\cal F}_\Gamma)\longrightarrow H^{0,\,1}_{[\Gamma],\,def}(X,\,T^{1,\,0}X)$ is {\it injective} if ${\cal F}_\Gamma$ and ${\cal G}_\Gamma$ are locally free and induce a direct-sum splitting $T^{1,\,0}X = {\cal F}_\Gamma\oplus{\cal G}_\Gamma$.

\vspace{1ex}

Our unobstructedness result is the following

\begin{The}\label{The:full-unobstructedness} Let $X$ be a compact complex manifold, with $\mbox{dim}_\C X = n = 2p+1$, equipped with a {\bf holomorphic $p$-contact structure} $\Gamma\in C^\infty_{p,\,0}(X,\,\C)$ such that ${\cal F}_\Gamma$ and ${\cal G}_\Gamma$ are {\bf locally free} and induce a direct-sum splitting $T^{1,\,0}X = {\cal F}_\Gamma\oplus{\cal G}_\Gamma$.

 Suppose that $(X,\,\Gamma)$ is a {\bf partially ${\cal F}_\Gamma$-directional page-$1$-$\partial\bar\partial$-manifold} and a {\bf partially vertically $\partial\bar\partial$-manifold}. Further suppose that the sheaf ${\cal F}_\Gamma$ is {\bf cohomologically integrable in bidegree $(0,\,1)$} and {\bf constantly maximally non-integrable in bidegree $(0,\,1)$}.

  Then, the {\bf essential horizontal deformations} of $X$ are {\bf unobstructed}.

\end{The}

\noindent {\it Proof.} It consists in a succession of steps that we spell out separately.

\vspace{1ex}

$\bullet$ {\it Step 1.} Let $[\theta_1]_{\bar\partial}\in E^{0,\,1}_2(X,\,{\cal F}_\Gamma)$ be arbitrary. Thanks to (1) of Proposition \ref{Prop:E_2_directional} and to the {\it partial ${\cal F}_\Gamma$-directional page-$1$-$\partial\bar\partial$-}assumption on $(X,\,\Gamma)$, there exists a {\it constantly horizontal} representative $\psi_1\in C^\infty_{0,\,1}(X,\,T^{1,\,0}X)$ of the class $[\theta_1]_{\bar\partial}$. This means that $\bar\partial\psi_1 = 0$, $[\psi_1]_{\bar\partial} = [\theta_1]_{\bar\partial}\in H^{0,\,1}_\Gamma(X,\,T^{1,\,0}X)$, $\psi_1\lrcorner\Gamma = 0$ and $\partial(\psi_1\lrcorner\partial\Gamma) = 0$.

Note that these properties imply $\partial(\psi_1\lrcorner u_\Gamma) = 0$. Indeed, $\psi_1\lrcorner u_\Gamma = (\psi_1\lrcorner\Gamma)\wedge\partial\Gamma + \Gamma\wedge(\psi_1\lrcorner\partial\Gamma) = \Gamma\wedge(\psi_1\lrcorner\partial\Gamma)$, hence \begin{eqnarray*}\partial(\psi_1\lrcorner u_\Gamma) = \partial\Gamma\wedge(\psi_1\lrcorner\partial\Gamma).\end{eqnarray*} This last product is zero since $\partial\Gamma\wedge\partial\Gamma = 0$ (for bidegree reasons), hence \begin{eqnarray*}0 = \psi_1\lrcorner(\partial\Gamma\wedge\partial\Gamma) = (\psi_1\lrcorner\partial\Gamma)\wedge\partial\Gamma + \partial\Gamma\wedge(\psi_1\lrcorner\partial\Gamma) = 2\,\partial\Gamma\wedge(\psi_1\lrcorner\partial\Gamma).\end{eqnarray*}

\vspace{1ex}

$\bullet$ {\it Step 2.} The vector-valued form $[\psi_1,\,\psi_1]\in C^\infty_{0,\,2}(X,\,T^{1,\,0}X)$ is {\it vertical} (i.e. it has the property $[\psi_1,\,\psi_1]\lrcorner\partial\Gamma = 0$), thanks to the constantly maximal non-integrability assumption in bidegree $(0,\,1)$ made on ${\cal F}_\Gamma$. Moreover, $[\psi_1,\,\psi_1]$ is also {\it $\bar\partial$-exact}. Indeed, its verticality implies \begin{eqnarray*}[\psi_1,\,\psi_1]\lrcorner u_\Gamma = \bigg([\psi_1,\,\psi_1]\lrcorner\Gamma\bigg)\wedge\partial\Gamma,\end{eqnarray*} which further implies that $[\psi_1,\,\psi_1]\lrcorner u_\Gamma$ is $\bar\partial$-exact since $[\psi_1,\,\psi_1]\lrcorner\Gamma$ is $\bar\partial$-exact (thanks to the cohomological integrability assumption in bidegree $(0,\,1)$ made on ${\cal F}_\Gamma$ and to $\psi_1\lrcorner\Gamma = 0$) and $\partial\Gamma$ is $\bar\partial$-closed. Now, the cohomological version (\ref{eqn:C_Y-isomorphism_cohom}) of the Calabi-Yau isomorphism $T_{[\Gamma]} : H^{0,\,2}_{\bar\partial}(X,\,T^{1,\,0}X)\longrightarrow H^{n-1,\,2}_{\bar\partial}(X,\,\C)$ implies that the $\bar\partial$-exactness of $[\psi_1,\,\psi_1]\lrcorner u_\Gamma$ is equivalent to the $\bar\partial$-exactness of $[\psi_1,\,\psi_1]$ (since $T_{[\Gamma]}([[\psi_1,\,\psi_1]]_{\bar\partial}) = 0$ if and only if $[[\psi_1,\,\psi_1]]_{\bar\partial}) = 0$).


We now notice that $[\psi_1,\,\psi_1]$ is even {\it constantly vertical}, namely it has the extra property \begin{eqnarray}\label{eqn:psi_1-bracket_constantly-vertical}\partial\bigg([\psi_1,\,\psi_1]\lrcorner\Gamma\bigg) = 0.\end{eqnarray} Indeed, the generalised Tian-Todorov-type formula (\ref{eqn:generalised_basic-trick_non-C-Y}) (in Lemma \ref{Lem:generalised_Tian-Todorov}) applied with $\alpha = \Gamma$ and $\theta_1 = \theta_2 = \psi_1$ yields: \begin{equation*}[\psi_1,\, \psi_1]\lrcorner\Gamma = -\partial\bigg(\psi_1\lrcorner(\psi_1\lrcorner\Gamma)\bigg) + 2\,\psi_1\lrcorner\partial(\psi_1\lrcorner\Gamma) - \psi_1\lrcorner(\psi_1\lrcorner\partial\Gamma)\end{equation*} after using the equality $L_{\psi_1}(\Gamma) = \partial(\psi_1\lrcorner\Gamma) - \psi_1\lrcorner\partial\Gamma$. Now, $\psi_1\lrcorner\Gamma = 0$ and $\psi_1\lrcorner(\psi_1\lrcorner\partial\Gamma)\in\ker\partial$ (see below), so the above equality implies $[\psi_1,\, \psi_1]\lrcorner\Gamma\in\ker\partial$, as stated in (\ref{eqn:psi_1-bracket_constantly-vertical}).

To see that $\psi_1\lrcorner(\psi_1\lrcorner\partial\Gamma)\in\ker\partial$, one writes: \begin{equation*}\partial\bigg(\psi_1\lrcorner(\psi_1\lrcorner\partial\Gamma)\bigg) = L_{\psi_1}(\psi_1\lrcorner\partial\Gamma) + \psi_1\lrcorner\partial(\psi_1\lrcorner\partial\Gamma) = \psi_1\lrcorner L_{\psi_1}(\partial\Gamma) - [\psi_1,\,\psi_1]\lrcorner\partial\Gamma,\end{equation*} where to get the latter equality we used the property $\partial(\psi_1\lrcorner\partial\Gamma) = 0$ (recall that $\psi$ is constantly horizontal) and the general commutation identity (\ref{eqn:Lie-derivatives-mixed_prop_2}) of Lemma \ref{Lem:Lie-derivatives-mixed_prop}. Now, $ [\psi_1,\,\psi_1]\lrcorner\partial\Gamma = 0$ (as seen above) and $L_{\psi_1}(\partial\Gamma) = \partial(\psi_1\lrcorner\partial\Gamma) - \psi_1\lrcorner\partial^2\Gamma = 0$, so we conclude that $\psi_1\lrcorner(\psi_1\lrcorner\partial\Gamma)\in\ker\partial$, as desired.

\vspace{1ex}

Since $[\psi_1,\,\psi_1]\in\mbox{Im}\,\bar\partial$ and $[\psi_1,\,\psi_1]\lrcorner\partial\Gamma = 0$, (iii)(b) of Proposition \ref{Prop:F-G-sheaves_properties} ensures the existence of $\psi'_2\in C^\infty_{0,\,1}(X,\,T^{1,\,0}X)$ such that \begin{eqnarray*}\bar\partial\psi'_2 = \frac{1}{2}\,[\psi_1,\,\psi_1] \hspace{3ex}\mbox{and}\hspace{3ex} \psi'_2\lrcorner\partial\Gamma = 0.\end{eqnarray*} In particular, we get \begin{eqnarray*}\bar\partial(\psi'_2\lrcorner\Gamma) = \frac{1}{2}\,[\psi_1,\,\psi_1]\lrcorner\Gamma, \hspace{3ex}\mbox{hence}\hspace{3ex} \bar\partial\partial(\psi'_2\lrcorner\Gamma) = - \frac{1}{2}\,\partial\bigg([\psi_1,\,\psi_1]\lrcorner\Gamma\bigg) = 0,\end{eqnarray*} where (\ref{eqn:psi_1-bracket_constantly-vertical}) was used to get the last equality.

Thus, $\psi'_2\in C^\infty_{0,\,1}(X,\,{\cal G}_\Gamma)$ has the property $\partial(\psi'_2\lrcorner\Gamma)\in\ker\bar\partial$. Consequently, the {\it partial vertical $\partial\bar\partial$-}hypothesis (see $(3)$(a) of Definition \ref{Def:dd-bar_F-G_directional}) made on $(X,\,\Gamma)$ ensures the existence of $\xi_2\in C^\infty(X,\,T^{1,\,0}X)$ such that \begin{eqnarray*}\xi_2\lrcorner\partial\Gamma = 0  \hspace{3ex}\mbox{and}\hspace{3ex} \partial\bar\partial(\xi_2\lrcorner\Gamma) = \partial(\psi'_2\lrcorner\Gamma).\end{eqnarray*}

Summing up, if we set $\psi_2:=\psi'_2 - \bar\partial\xi_2\in C^\infty_{0,\,1}(X,\,T^{1,\,0}X)$, we have \begin{eqnarray*}\psi_2\lrcorner\partial\Gamma = 0 \hspace{3ex}\mbox{and}\hspace{3ex} \partial(\psi_2\lrcorner\Gamma) = 0\end{eqnarray*} (two conditions that together mean that $\psi_2$ is {\it constantly vertical}), as well as \begin{eqnarray*}\label{eqn:integrability-condition_2-bis}\bar\partial\psi_2 = \frac{1}{2}\,[\psi_1,\,\psi_1]   \hspace{6ex} (\mbox{Eq.}\,\,(2)).\end{eqnarray*}

In other words, we have proved the existence of a {\it constantly vertical} solution $\psi_2\in C^\infty_{0,\,1}(X,\,T^{1,\,0}X)$ of equation (\ref{eqn:integrability-condition_nu}) for $\nu=2$ (i.e. $(\mbox{Eq.}\,\,(2))$).

\vspace{1ex}

 We go on to notice the following general fact.

\begin{Lem}\label{Lem:del-exactness_constant-verticality} If $\psi\in C^\infty_{0,\,1}(X,\,T^{1,\,0}X)$ is {\bf constantly vertical}, we have \begin{eqnarray}\label{eqn:del-exactness_constant-verticality}\psi\lrcorner u_\Gamma = \partial\bigg(-(\psi\lrcorner\Gamma)\wedge\Gamma\bigg).\end{eqnarray}

\end{Lem}

\noindent {\it Proof.} The constant verticality assumption on $\psi$ means that $\psi\lrcorner\partial\Gamma = 0$ and $\partial(\psi\lrcorner\Gamma) = 0$. Since $u_\Gamma = \Gamma\wedge\partial\Gamma$, we get: \begin{eqnarray*}\psi\lrcorner u_\Gamma = (\psi\lrcorner\Gamma)\wedge\partial\Gamma + \Gamma\wedge(\psi\lrcorner\partial\Gamma) =  (\psi\lrcorner\Gamma)\wedge\partial\Gamma = \partial\bigg(-(\psi\lrcorner\Gamma)\wedge\Gamma\bigg).\end{eqnarray*} This proves the contention.  \hfill $\Box$

\vspace{1ex}

In particular, in our case we get $\psi_2\lrcorner u_\Gamma\in\mbox{Im}\,\partial$, hence $\psi_2\lrcorner u_\Gamma\in\ker\partial$.

\vspace{1ex}

$\bullet$ {\it Step 3.} We need to prove the existence of a solution $\psi_3\in C^\infty_{0,\,1}(X,\,T^{1,\,0}X)$ of the equation \begin{eqnarray*}\label{eqn:integrability-condition_3}\bar\partial\psi_3 = [\psi_1,\,\psi_2]   \hspace{6ex} (\mbox{Eq.}\,\,(3)).\end{eqnarray*} (See (\ref{eqn:integrability-condition_nu}) for $\nu = 3$.) However, in our current setting, $\psi_1$ is {\it constantly horizontal} and $\psi_2$ is {\it constantly vertical}, so $[\psi_1,\,\psi_2] = 0$ thanks to (i) of Corollary \ref{Cor:brackets_horizontal-vertical_forms}. Thus, we can choose $\psi_3 = 0$.

\vspace{1ex}

$\bullet$ {\it Step 4.} We now prove by induction on the integer $\nu\geq 2$ that there is a sequence $(\psi_\nu)_{\nu\geq 2}$ of solutions to the inductively defined sequence (\ref{eqn:integrability-condition_nu}) of equations such that: \begin{eqnarray}\label{eqn:induction-result}\nonumber & & \psi_{2k+1} = 0   \hspace{3ex} \mbox{for every} \hspace{1ex} k\geq 1; \\
  & & \psi_{2k} \hspace{1ex} \mbox{is {\bf constantly vertical} and satisfies} \hspace{1ex} \psi_{2k}\lrcorner u_\Gamma\in\mbox{Im}\,\partial \hspace{3ex} \mbox{for every} \hspace{1ex} k\geq 1.\end{eqnarray}

\vspace{1ex}

The cases $\nu = 2$ and $\nu = 3$ have been proved under Steps $2$ and $3$ above.

\vspace{1ex}

(a)\, If $\nu$ is odd (let then $\nu = 2l+1$), we need to show that $\psi_\nu = 0$ is a solution of the equation \begin{eqnarray*}\bar\partial\psi_\nu = [\psi_1,\,\psi_{\nu-1}] + \dots + [\psi_l,\,\psi_{\nu-l}]   \hspace{6ex} (\mbox{Eq.}\,\,(\nu)).\end{eqnarray*} (This is what (\ref{eqn:integrability-condition_nu}) becomes after using the symmetries $[\psi_s,\,\psi_{\nu-s}] = [\psi_{\nu-s},\,\psi_s]$ for all $s$.) This is so because the right-hand side of this equation vanishes. Indeed, by (i) of Corollary \ref{Cor:brackets_horizontal-vertical_forms}, $[\psi_1,\,\psi_{\nu-1}] = 0$ since $\psi_1$ is constantly horizontal and, by the induction hypothesis, $\psi_{\nu-1} = \psi_{2l}$ is constantly vertical. Meanwhile, $[\psi_s,\,\psi_{\nu-s}] = 0$ for every $s=2,\dots, l$ since one of the indices $s$ and $\nu-s$ is odd and, therefore, either $\psi_s = 0$ or $\psi_{\nu-s} = 0$ by the induction hypothesis.

\vspace{1ex}

(b)\, Suppose that $\nu$ is even. Let $\nu = 2l$. We need to find a solution $\psi_\nu$ to the equation \begin{eqnarray*}\bar\partial\psi_\nu = [\psi_1,\,\psi_{\nu-1}] + \dots + [\psi_{l-1},\,\psi_{l+1}] + \frac{1}{2}\,[\psi_l,\,\psi_l]  \hspace{6ex} (\mbox{Eq.}\,\,(\nu)).\end{eqnarray*} (This is what (\ref{eqn:integrability-condition_nu}) becomes after using the symmetries $[\psi_s,\,\psi_{\nu-s}] = [\psi_{\nu-s},\,\psi_s]$ for all $s$.) Since $\nu- 1$, $\nu-3$, etc are odd, the induction hypothesis ensures that $[\psi_1,\,\psi_{\nu-1}] = [\psi_3,\,\psi_{\nu-3}] = \dots = 0$. Thus, the right-hand side of $(\mbox{Eq.}\,\,(\nu))$ reduces to \begin{eqnarray*}R_\nu & := & \sum\limits_{s=1}^{r-1}[\psi_{2s},\,\psi_{\nu-2s}] + \frac{1}{2}\,[\psi_{2r},\,\psi_{2r}],    \hspace{6ex} \mbox{if}\hspace{1ex} l=2r \\
  R_\nu & := & \sum\limits_{s=1}^r[\psi_{2s},\,\psi_{\nu-2s}],    \hspace{6ex} \mbox{if}\hspace{1ex} l=2r+1.\end{eqnarray*}

Fix any $s$. Since $\psi_{2s}$ and $\psi_{\nu-2s}$ are both vertical (by the induction hypothesis) and since ${\cal G}_\Gamma$ is integrable (see Observation \ref{Obs:G_Gamma_integrability}), $[\psi_{2s},\,\psi_{\nu-2s}]$ is vertical. This implies that $R_\nu\lrcorner\partial\Gamma = 0$.

Meanwhile, the generalised Tian-Todorov formula (\ref{eqn:generalised_basic-trick_non-C-Y}) (of (ii) in Lemma \ref{Lem:generalised_Tian-Todorov}) applied with $\alpha = \Gamma$ and $\theta_1 = \psi_{2s}$, $\theta_2 = \psi_{\nu-2s} = \psi_{2(l-s)}$ yields: \begin{equation*}[\psi_{2s},\, \psi_{2(l-s)}]\lrcorner\Gamma = -\partial\bigg(\psi_{2s}\lrcorner(\psi_{2(l-s)}\lrcorner\Gamma)\bigg) + \psi_{2s}\lrcorner\partial(\psi_{2(l-s)}\lrcorner\Gamma)  + \psi_{2(l-s)}\lrcorner\partial(\psi_{2s}\lrcorner\Gamma) - \psi_{2s}\lrcorner(\psi_{2(l-s)}\lrcorner\partial\Gamma)\end{equation*} after using the equality $L_{\psi_{2s}}(\Gamma) = \partial(\psi_{2s}\lrcorner\Gamma) - \psi_{2s}\lrcorner\partial\Gamma$ and its analogue for $\psi_{2(l-s)}$. Since $\psi_{2s}$ and $\psi_{2(l-s)}$ are constantly vertical (by the induction hypothesis), $\partial(\psi_{2s}\lrcorner\Gamma) = \psi_{2(l-s)}\lrcorner\Gamma = 0$. Moreover, $\psi_{2(l-s)}\lrcorner\partial\Gamma = 0$ since $\psi_{2(l-s)}$ is vertical.

On the one hand, we conclude that $[\psi_{2s},\, \psi_{2(l-s)}]\lrcorner\Gamma\in\mbox{Im}\,\partial$, hence $[\psi_{2s},\, \psi_{2(l-s)}]\lrcorner\Gamma\in\ker\partial$. This, together with the verticality of $[\psi_{2s},\, \psi_{\nu-2s}]$ (seen above), means that $[\psi_{2s},\,\psi_{\nu-2s}]$ is constantly vertical.

On the other hand, for the right-hand side $R_\nu$ of $(\mbox{Eq.}\,\,(\nu))$ we conclude that \begin{eqnarray*}R_\nu\lrcorner\Gamma & := & -\partial\bigg(\sum\limits_{s=1}^{r-1}\psi_{2s}\lrcorner(\psi_{2(l-s)}\lrcorner\Gamma) - \frac{1}{2}\,\psi_{2r}\lrcorner(\psi_{2r}\lrcorner\Gamma)\bigg), \hspace{6ex} \mbox{if}\hspace{1ex} l=2r \\
  R_\nu\lrcorner\Gamma & := & -\partial\bigg(\sum\limits_{s=1}^r\psi_{2s}\lrcorner(\psi_{2(l-s)}\lrcorner\Gamma)\bigg), \hspace{6ex} \mbox{if}\hspace{1ex} l=2r+1.\end{eqnarray*}

Meanwhile, $R_\nu\in\ker\bar\partial$ (a general well-known fact in all these equations). Since $\bar\partial\Gamma = 0$, this yields $R_\nu\lrcorner\Gamma\in\ker\bar\partial$. Together with the above shape of $R_\nu\lrcorner\Gamma$ as $\partial(\dots)$, this implies, thanks to the {\it partial vertical $\partial\bar\partial$-}hypothesis (see $(3)$(b) of Definition \ref{Def:dd-bar_F-G_directional}) made on $(X,\,\Gamma)$, that $R_\nu\lrcorner\Gamma\in\mbox{Im}\,\bar\partial$.

Therefore, using also the verticality of $R_\nu$ (seen above), we get $R_\nu\lrcorner u_\Gamma = (R_\nu\lrcorner\Gamma)\wedge\partial\Gamma\in\mbox{Im}\,\bar\partial$. Thanks to the cohomological version (\ref{eqn:C_Y-isomorphism_cohom}) of the Calabi-Yau isomorphism, this is equivalent to $R_\nu\in\mbox{Im}\,\bar\partial$, which means that equation $(\mbox{Eq.}\,\,(\nu))$ is solvable.

The existence of a solution $\psi_\nu$ that is constantly vertical and satisfies $\psi_\nu\lrcorner u_\Gamma\in\mbox{Im}\,\partial$ is proved in the same way as in the case of $\psi_2$ in the last part of Step $2$:

\vspace{1ex}

-we first use (iii)(b) of Proposition \ref{Prop:F-G-sheaves_properties} to get a vertical solution $\psi'_\nu$ to equation $(\mbox{Eq.}\,\,(\nu))$;

\vspace{1ex}

-we then use the {\it partial vertical $\partial\bar\partial$-}hypothesis (see $(3)$(a) of Definition \ref{Def:dd-bar_F-G_directional}) made on $(X,\,\Gamma)$ to correct $\psi'_\nu$ to a constantly vertical solution $\psi_\nu:=\psi'_\nu - \bar\partial\xi_\nu\in C^\infty_{0,\,1}(X,\,T^{1,\,0}X)$ to equation $(\mbox{Eq.}\,\,(\nu))$ for some vertical $\xi_\nu\in C^\infty(X,\,T^{1,\,0}X)$ such that $\partial\bar\partial(\xi_\nu\lrcorner\Gamma) = \partial(\psi'_\nu\lrcorner\Gamma)$;

\vspace{1ex}

-finally, we use the general Lemma \ref{Lem:del-exactness_constant-verticality}
to get $\psi_\nu\lrcorner u_\Gamma\in\mbox{Im}\,\partial$, hence $\psi_\nu\lrcorner u_\Gamma\in\ker\partial$.

\vspace{2ex}

This completes the induction process and the proof of Theorem \ref{The:full-unobstructedness}.  \hfill $\Box$

\vspace{3ex}

\noindent {\bf Acknowledgements.} The authors are grateful to F. Campana for pointing out to one of them the reference [EMS77] and for an enlightening exchange of emails, as well as to S. Boucksom for useful discussions. The second-named author visited his co-authors separately in their respective institutions and is grateful to the Osaka University and to the University of Zaragoza for support and hospitality. The third-named author  was partially supported by grant PID2023-148446NB-I00, funded by MICIU/AEI/10.13039/501100011033, 
and by grant E22-23R ``Algebra y Geometr\' ia'' (Gobierno de Arag\'on/FEDER).

\vspace{3ex}

\noindent {\bf References.} \\



\noindent [BW58]\, W.M. Boothby, H.C. Wang --- {\it On Contact Manifolds} --- Ann. of Math. {\bf 68} (1958) 721--734.

\vspace{1ex}

\noindent [CP94]\, F. Campana, T. Peternell --- {\it Cycle spaces} --- Several Complex Variables, VII,  319-349, Encyclopaedia Math. Sci., {\bf 74}, Springer, Berlin (1994). 

\vspace{1ex}

\noindent [CFGU97]\, L.A. Cordero, M. Fern\'andez, A.Gray, L. Ugarte --- {\it A General Description of the Terms in the Fr\"olicher Spectral Sequence} --- Differential Geom. Appl. {\bf 7} (1997), no. 1, 75--84.

\vspace{1ex}

\noindent [Dem02]\, J.-P. Demailly --- {\it On the Frobenius Integrability of Certain Holomorphic $p$-Forms } --- In: Bauer, I., Catanese, F., Peternell, T., Kawamata, Y., Siu, YT. (eds) Complex Geometry. Springer, Berlin, Heidelberg. \url{https://doi.org/10.1007/978-3-642-56202-0_6}

\vspace{1ex}

\noindent [Dem97]\, J.-P. Demailly --- {\it Complex Analytic and Algebraic Geometry} --- http://www-fourier.ujf-grenoble.fr/~demailly/books.html

\vspace{1ex}

\noindent [EMS77]\, R. Edwards, K. Millett, D. Sullivan --- {\it Foliations with All Leaves Compact} --- Topology {\bf 16} (1977), 13-32.





\vspace{1ex}

\noindent [For00]\, B. Foreman --- {\it Boothby-Wang Fibrations on Complex Contact Manifolds} --- Differential Geometry and its Applications {\bf 13} (2000), 179-196. 



\vspace{1ex}


\noindent [KPU25]\, H. Kasuya, D. Popovici, L. Ugarte --- {\it Properties of Holomorphic $p$-Contact Manifolds} ---  arXiv:2511.10818v1 [math.DG]. 

\vspace{1ex}

\noindent [Kob87]\, S. Kobayashi --- {\it Differential Geometry of Complex Vector Bundles} --- Princeton University Press, 1987.

\vspace{1ex}

\noindent [Kod86]\, K. Kodaira --- {\it Complex Manifolds and Deformations of Complex Structures} --- Grundlehren der Math. Wiss. {\bf 283}, Springer (1986).

\vspace{1ex}

\noindent [Lie78]\, D. Lieberman --- {\it Compactness of the Chow Scheme: Applications to Automorphisms and Deformations of K\"ahler Manifolds} --- Lect. Notes Math. {\bf 670} (1978), 140-186.















\vspace{1ex}

\noindent [Pal57]\, R.S. Palais --- {\it A Global Formulation of the Lie Theory of Transformation Groups} --- Mem. Amer. Math. Soc. {\bf 22} (Amer. Math. Soc., Providence, 1957).

\vspace{1ex}

\noindent [Pop18]\, D. Popovici --- {\it Non-K\"ahler Mirror Symmetry of the Iwasawa Manifold} --- Int. Math. Res. Not. IMRN 2020, no. {\bf 23}, 9471--9538, doi:10.1093/imrn/rny256.

\vspace{1ex}

\noindent [Pop19]\,  D. Popovici --- {\it Holomorphic Deformations of Balanced Calabi-Yau $\partial\bar\partial$-Manifolds} --- Ann. Inst. Fourier, {\bf 69} (2019) no. 2, pp. 673-728. doi : 10.5802/aif.3254.

\vspace{1ex}

\noindent [PSU20a]\,  D. Popovici, J. Stelzig, L. Ugarte --- {\it Higher-Page Hodge Theory of Compact Complex Manifolds} --- arXiv e-print AG 2001.02313v3; Ann. Sc. Norm. Super. Pisa Cl. Sci. (5) Vol. XXV (2024), 1431-1464. 

\vspace{1ex}

\noindent [PSU20b]\, D. Popovici, J. Stelzig, L. Ugarte --- {\it Higher-Page Bott-Chern and Aeppli Cohomologies and Applications} --- J. reine angew. Math., doi: 10.1515/crelle-2021-0014.

\vspace{1ex}

\noindent [PSU20c]\, D. Popovici, J. Stelzig, L. Ugarte ---  {\it Deformations of Higher-Page Analogues of $\partial\bar\partial$-Manifolds} --- Math. Z. (2021), https://doi.org/10.1007/s00209-021-02861-0.

\vspace{1ex}

\noindent [PU23]\,  D. Popovici, L. Ugarte --- {\it A Moment Map for the Space of Maps to a Balanced Manifold} --- arXiv e-print DG 2311.00485v1



\vspace{1ex}

\noindent [Tia87]\, G. Tian --- {\it Smoothness of the Universal Deformation Space of Compact Calabi-Yau Manifolds and Its Petersson-Weil Metric} --- Mathematical Aspects of String Theory (San Diego, 1986), Adv. Ser. Math. Phys. 1, World Sci. Publishing, Singapore (1987), 629--646.

\vspace{1ex}

\noindent [Tod89]\, A. N. Todorov --- {\it The Weil-Petersson Geometry of the Moduli Space of $SU(n\geq 3)$ (Calabi-Yau) Manifolds I} --- Comm. Math. Phys. {\bf 126} (1989), 325-346.



\vspace{3ex}

\noindent Graduate School of Mathematics    \hfill Institut de Math\'ematiques de Toulouse

\noindent Nagoya University, Furocho, Chikusaku,    \hfill  Universit\'e de Toulouse

\noindent Nagoya, Japan, 464-8602  \hfill  118 route de Narbonne, 31062 Toulouse, France

\noindent  Email: kasuya@math.nagoya-u.ac.jp                \hfill     Email: popovici@math.univ-toulouse.fr

\vspace{2ex}

\noindent and

\vspace{2ex}

\noindent Departamento de Matem\'aticas\,-\,I.U.M.A., Universidad de Zaragoza,

\noindent Campus Plaza San Francisco, 50009 Zaragoza, Spain

\noindent Email: ugarte@unizar.es

\end{document}